\documentclass[11pt,twoside]{article}

\usepackage{fullpage}

\usepackage{epsf}
\usepackage{fancyhdr}
\usepackage{graphics}
\usepackage{graphicx}
\usepackage{psfrag}

\usepackage[linesnumbered,ruled]{algorithm2e}% http://ctan.org/pkg/algorithm2e
\DontPrintSemicolon

\usepackage{color}

\usepackage{amsthm}
\usepackage{amsfonts}
\usepackage{amsmath}
\usepackage{amssymb,bbm}

\usepackage[colorlinks]{hyperref}
\usepackage{nicefrac}

\usepackage{chngpage}

 \usepackage{tabularx}%

\usepackage{enumitem}
\usepackage{booktabs}
\usepackage{caption,subcaption}

\usepackage{tikz}
\usetikzlibrary{positioning}
\usetikzlibrary{arrows.meta}

\usepackage{mathtools}

\newcommand{\denb}{\varrho}
\newcommand{\ttensor}{\mathcal{T}}
\newcommand{\spec}[1]{\rho\parenth{#1}}
\newcommand{\isoperi}{\psi}

\newcommand{\co}{A}
\newcommand{\qo}{Q}

\newcommand{\potentwo}{\Gamma}

\newcommand{\indicator}{\mathbbm{1}}

\newcommand{\defn}{:=}

\newtheoremstyle{named}{}{}{\itshape}{}{\bfseries}{.}{.5em}{\thmnote{#3's }#1}
\theoremstyle{named}

\theoremstyle{plain}

\newtheorem{theorem}{Theorem}

\newtheorem{lemma}{Lemma}

\newtheorem{conjecture}{Conjecture}

\newlength{\widebarargwidth}
\newlength{\widebarargheight}
\newlength{\widebarargdepth}

\makeatletter
\long\def\@makecaption#1#2{
        \vskip 0.8ex
        \setbox\@tempboxa\hbox{\small {\bf #1:} #2}
        \parindent 1.5em  %% How can we use the global value of this???
        \dimen0=\hsize
        \advance\dimen0 by -3em
        \ifdim \wd\@tempboxa >\dimen0
                \hbox to \hsize{
                        \parindent 0em
                        \hfil
                        \parbox{\dimen0}{\def\baselinestretch{0.96}\small
                                {\bf #1.} #2
                                }
                        \hfil}
        \else \hbox to \hsize{\hfil \box\@tempboxa \hfil}
        \fi
        }
\makeatother

\long\def\comment#1{}
\definecolor{battleshipgrey}{rgb}{0.52, 0.52, 0.51}
\definecolor{darkgray}{rgb}{0.66, 0.66, 0.66}
\definecolor{darkgreen}{rgb}{0.0, 0.2, 0.13}
\definecolor{darkspringgreen}{rgb}{0.09, 0.45, 0.27}
\definecolor{dukeblue}{rgb}{0.0, 0.0, 0.61}
\definecolor{olivedrab7}{rgb}{0.24, 0.2, 0.12}
\definecolor{darkblue}{rgb}{0.0, 0.0, 0.55}
\definecolor{darkscarlet}{rgb}{0.34, 0.01, 0.1}
\definecolor{candyapplered}{rgb}{1.0, 0.03, 0.0}
\definecolor{ao(english)}{rgb}{0.0, 0.5, 0.0}
\definecolor{applegreen}{rgb}{0.55, 0.71, 0.0}

\newcommand{\vecnorm}[2]{\left\| #1\right\|_{#2}}

\newcommand{\Exs}{\ensuremath{{\mathbb{E}}}}
\newcommand{\Prob}{\ensuremath{{\mathbb{P}}}}

\newcommand{\usedim}{\ensuremath{d}}
\newcommand{\dims}{\usedim}

\DeclareMathOperator{\Var}{Var}
\DeclareMathOperator{\Cov}{Cov}
\DeclareMathOperator{\trace}{Tr}

\newcommand{\Normal}{\ensuremath{\mathcal{N}}}

\newcommand{\Ind}{\ensuremath{\mathbb{I}}}

\newcommand{\real}{\ensuremath{\mathbb{R}}}

\newcommand{\interior}[1]{%
  {\kern0pt#1}^{\mathrm{o}}%
}

\newcommand{\brackets}[1]{\left[ #1 \right]}
\newcommand{\parenth}[1]{\left( #1 \right)}

\newcommand{\braces}[1]{\left\{ #1 \right \}}
\newcommand{\abss}[1]{\left| #1 \right |}

\newcommand{\tp}{^\top}

\begin{document}

%%%%%%% TITLE PAGE %%%%%%%%%%%%%%%%%%%%%%%%%%%%%%%%%%%%%%%%%%%%%%%%%%%

\begin{center}

{\bf{\LARGE{An Almost Constant Lower Bound of the Isoperimetric Coefficient in the KLS Conjecture}}}

\vspace*{.2in}

{\large{
\begin{tabular}{cccc}
Yuansi Chen\\
\end{tabular}
}}

\vspace*{.2in}

\begin{tabular}{c}
Seminar for Statistics \\
ETH, Z\"urich
\end{tabular}

\vspace*{.2in}

\today

\vspace*{.2in}

\begin{abstract}
We prove an almost constant lower bound of the isoperimetric coefficient in the KLS conjecture. The lower bound has the dimension dependency $\dims^{-o_d(1)}$. When the dimension is large enough, our lower bound is tighter than the previous best bound which has the dimension dependency $\dims^{-1/4}$. Improving the current best lower bound of the isoperimetric coefficient in the KLS conjecture has many implications, including improvements of the current best bounds in Bourgain's slicing conjecture and in the thin-shell conjecture, better concentration inequalities for Lipschitz functions of log-concave measures and better mixing time bounds for MCMC sampling algorithms on log-concave measures.
\end{abstract}

\end{center}

%!TEX root = main_paper.tex

\section{Introduction} % (fold)
\label{sec:introduction}
Given a distribution, the isoperimetric coefficient of a subset is the ratio of the measure of the subset boundary to the minimum of the measures of the subset and its complement. Taking the minimum of such ratios over all subsets defines the isoperimetric coefficient of the distribution, also called the Cheeger isoperimetric coefficient of the distribution.

Kannan, Lov\'asz and Simonovits (KLS)~\cite{kannan1995isoperimetric} conjecture that for any distribution that is log-concave, the Cheeger isoperimetric coefficient equals to that achieved by half-spaces up to a universal constant factor. If the conjecture is true, the Cheeger isoperimetric coefficient can be determined by going through all the half-spaces instead of all subsets. For this reason, the KLS conjecture is also called the KLS hyperplane conjecture. To make it precise, we start by formally defining log-concave distributions and then we state the conjecture.

A probability density function $p: \real^\dims \rightarrow \real$ is \textit{log-concave} if its logarithm is concave, i.e., for any $x, y \in \real^{\dims} \times \real^{\dims}$ and for any $\lambda \in [0, 1]$,
\begin{align}
  \label{eq:def_logconcave}
  p(\lambda x + (1 - \lambda) y) \geq p(x)^\lambda p(y)^{1-\lambda}.
\end{align}
Common probability distributions such as Gaussian, exponential and logistic are log-concave. This definition also includes any uniform distribution over a \textit{convex set} defined as follows. A subset $K \subset \real^\dims$ is \textit{convex} if $\forall x, y \in K \times K, z \in [x, y] \implies z \in K$.
The \textit{isoperimetric coefficient} $\isoperi(p)$ of a density $p$ in $\real^\dims$ is defined as
\begin{align}
  \label{eq:def_isoperi}
  \isoperi(p) \defn \inf_{S \subset \real^\dims }\frac{p^+(\partial S)}{\min(p(S), p(S^c))}
\end{align}
where $p(S) = \int_{x \in S} p(x) dx$ and the boundary measure of the subset is
\begin{align*}
  p^+(\partial S) \defn \underset{\epsilon \rightarrow 0^+}{\lim\inf}\ \frac{p\parenth{\braces{x: \mathbf{d}(x, S) \leq \epsilon}} - p(S)}{\epsilon},
\end{align*}
where $\mathbf{d}(x, S)$ is the Euclidean distance between $x$ and the subset $S$.

The KLS conjecture is stated by Kannan, Lov\'asz and Simonovits~\cite{kannan1995isoperimetric} as follows.
\begin{conjecture}
  \label{cj:main}
  There exists a universal constant $c$, such that for any log-concave density $p$ in $\real^\dims$, we have
  \begin{align*}
    \isoperi(p) \geq \frac{c}{\sqrt{\spec{p}}},
  \end{align*}
  where $\spec{p}$ is the spectral norm of the covariance matrix of $p$. In other words, $\spec{p} = \vecnorm{A}{2}$, where $A = \Cov_{X \sim p} (X)$ is the covariance matrix.
\end{conjecture}
An upper bound of $\isoperi(p)$ of the same form is relatively easy and it was shown to be achieved by half-spaces~\cite{kannan1995isoperimetric}. Proving the lower bound on $\isoperi(p)$ up to some small factors in Conjecture~\ref{cj:main} is the main goal of this paper. We say a log-concave density is \textit{isotropic} if its mean $\Exs_{X\sim p} [X]$ equals to $0$ and its covariance $\Cov_{X\sim p}(X)$ equals to $\Ind_\dims$. In the case of isotropic log-concave densities, the KLS conjecture states that any isotropic log-concave density has its isoperimetric coefficient lower bounded by a universal constant.

There are many attempts trying to lower bound the Cheeger isoperimetric coefficient in the KLS conjecture. We refer readers to the survey paper by Lee and Vempala~\cite{lee2018kannan} for a detailed exposition of these attempts. In particular, the original KLS paper~\cite{kannan1995isoperimetric} (Theorem 5.1) shows that for any log-concave density $p$ with covariance matrix $A$,
\begin{align}
  \label{eq:KLS_original_bound}
  \isoperi(p) \geq \frac{\log(2)}{\sqrt{\trace\parenth{A}}}.
\end{align}
The original KLS paper~\cite{kannan1995isoperimetric} only deals with uniform distributions over convex sets, but their proof techniques can be easily extended to show that the same results hold for all log-concave densities.
Remark that Equation~\eqref{eq:KLS_original_bound} implies $\isoperi(p) \geq \frac{\log(2)}{\dims^{1/2} \cdot \sqrt{\spec{p}}}$.
The current best bound is shown in Lee and Vempala~\cite{lee2016eldan}, where they show that there exists a universal constant $c$ such that for any log-concave density $p$ with covariance matrix $A$,
\begin{align}
  \label{eq:Lee_and_Vempala_bound}
  \isoperi(p) \geq \frac{c}{\parenth{\trace\parenth{A^2}}^{1/4}}.
\end{align}
It implies that $\isoperi(p) \geq \frac{c}{\dims^{1/4} \cdot \sqrt{\spec{p}}}$. Note that in Lee and Vempala~\cite{lee2016eldan}, their notation of $\isoperi(p)$ is the reciprocal of ours and it is later switched in Theorem 32 of the survey paper~\cite{lee2018kannan} by the same authors. As a result, the above bound is not a misstatement of the results in Lee and Vempala~\cite{lee2016eldan} and it is simply translated into our notations. In this paper, we improve the dimension dependency $\dims^{-1/4}$ to $\dims^{-o_\dims(1)}$ in the lower bound of the isoperimetric coefficient.

There are many implications of improving the lower bound in the KLS conjecture. The two closely related conjectures are Bourgain's slicing conjecture~\cite{bourgain1986high,ball1988logarithmically} and the thin-shell conjecture~\cite{anttila2003central}. It is worth noting that Bourgain~\cite{bourgain1986high} stated the slicing conjecture earlier than the introduction of the KLS conjecture. In terms of their connections to the KLS conjecture, Eldan and Klartag~\cite{eldan2011approximately} proved that the thin-shell conjecture implies Bourgain's slicing conjecture up to a universal constant factor. Later, Eldan~\cite{eldan2013thin} showed that the inverse of an lower bound of the isoperimetric coefficient is equivalent to an upper bound of the thin-shell constant in the thin-shell conjecture. Combining these two results, we have that an lower bound in the KLS conjecture implies upper bounds in the thin-shell conjecture and in Bourgain's slicing conjecture.

The current best upper bound of the thin-shell constant has the dimension dependency $\dims^{1/4}$ due to Lee and Vempala's~\cite{lee2016eldan} improvement in the KLS conjecture. The current best bound of the slicing constant in Bourgain's slicing conjecture also has the dimension dependency $\dims^{1/4}$, proved by Klartag~\cite{klartag2006convex} without using the KLS conjecture. Klartag's slicing constant bound is a slight improvement over Bourgain’s earlier slicing bound~\cite{bourgain1986high} which has the dimension dependency $\dims^{1/4}\log(\dims)$. Given the current best bounds in these three conjectures and the relation among them, we conclude that improving the current best lower bound in the KLS conjecture improves the current best bounds for the other two conjectures, as noted in Lee and Vempala~\cite{lee2018kannan}. For a detailed exposition of the three conjectures and related results since the introduction of Bourgain's slicing conjecture, we refer readers to Klartag and Milman~\cite{klartag2021slicing}.

Additionally, improving the lower bound in the KLS conjecture also improves concentration inequalities for Lipschitz functions of log-concave measures. It also leads to faster mixing time bounds of Markov chain Monte Carlo (MCMC) sampling algorithms on log-concave measures. Despite the great importance of these results, deriving these results from our new bound in the KLS conjecture is not the main focus of our paper. We refer readers to Milman~\cite{milman2009role} and Lee and Vempala~\cite{lee2018kannan} for more details about the abundant implications of the KLS conjecture.

\paragraph{Notation:} For two sequences $a_n$ and $b_n$ indexed by an integer $n$, we say that $a_n = o_n(b_n)$ if $\lim_{n \to \infty} \frac{a_n}{b_n} = 0$. The Euclidean norm of a vector $x \in \real^\dims$ is denoted by $\vecnorm{x}{2}$. The spectral norm of a square matrix $A \in \real^{\dims \times \dims}$ is denoted by $\vecnorm{A}{2}$. The Euclidean ball with center $x$ and radius $r$ is denoted by $\mathbb{B}(x, r)$. For a real number $x \in \real$, we denote its ceiling by $\lceil x \rceil = \min\braces{m \in \mathbb{Z} \mid m \geq x}$. We say a density $p$ is more log-concave than a Gaussian density $\varphi$ if $p$ can be written as a product form $p = \nu \cdot \varphi$ where $\varphi$ is the Gaussian density and $\nu$ is a log-concave function (that is, $\nu$ is proportional to a log-concave density). For a martingale $(M_t,\ t \in \real_+)$, we use $\brackets{M}_t$ to denote its \textit{quadratic variation}, defined as
\begin{align*}
  \brackets{M}_t = \sup_{k \in \mathbb{N}} \sup_{0 \leq t_1 \leq \cdots \leq t_k \leq t} \sum_{i=1}^k \parenth{M_{t_i} - M_{t_{i-1}}}^2.
\end{align*}

% section introduction (end)
%!TEX root = main_paper.tex

\section{Main results} % (fold)
\label{sec:main_results}
We prove the following lower bound on the isoperimetric coefficient of any log-concave density.
\begin{theorem}
  \label{thm:main_KLS}
  There exists a universal constant $c$ such that for any log-concave density $p$ in $\real^\dims$ and any integer $\ell \geq 1$, we have
  \begin{align}
    \label{eq:main_KLS}
    \isoperi(p) \geq \frac{1}{\brackets{c \cdot \ell \parenth{\log(\dims)+1}}^{\ell/2} \dims^{16/\ell} \cdot \sqrt{\spec{p}}}
  \end{align}
  where $\spec{p}$ is the spectral norm of the covariance matrix of $p$.
\end{theorem}
As a corollary, take $\ell = \left\lceil \parenth{\frac{\log(\dims)}{\log\log(\dims)}}^{1/2} \right\rceil$, then there exists a constant $c'$ such that
\begin{align*}
  \isoperi(p) \geq \frac{1}{\dims^{c' \parenth{\frac{\log\log(\dims)}{\log{\dims}}}^{1/2}}  \cdot \sqrt{\spec{p}}}.
\end{align*}
Since $\lim_{\dims \to \infty} \frac{\log\log(\dims)}{\log(\dims)} = 0$, for $\dims$ large enough, the above lower bound is better than any lower bound of the form $\frac{1}{\dims^{c''} \sqrt{\spec{p}}} $ ($c''$ is a positive constant) in terms of dimension $\dims$ dependency.

The proof of the main theorem uses the stochastic localization scheme introduced by Eldan~\cite{eldan2013thin}. Eldan uses this stochastic localization scheme to show that the thin shell conjecture is equivalent to the KLS conjecture up to a logarithmic factor. The construction of stochastic localization scheme uses elementary properties of semimartingales and stochastic integration. The main idea of Eldan's proof to derive the KLS conjecture from the thin shell conjecture is to smoothly multiply a Gaussian part to the log-concave density, so that the modified density is more log-concave than a Gaussian density. When the Gaussian part is large enough, one can then easily prove the isoperimetric inequality.

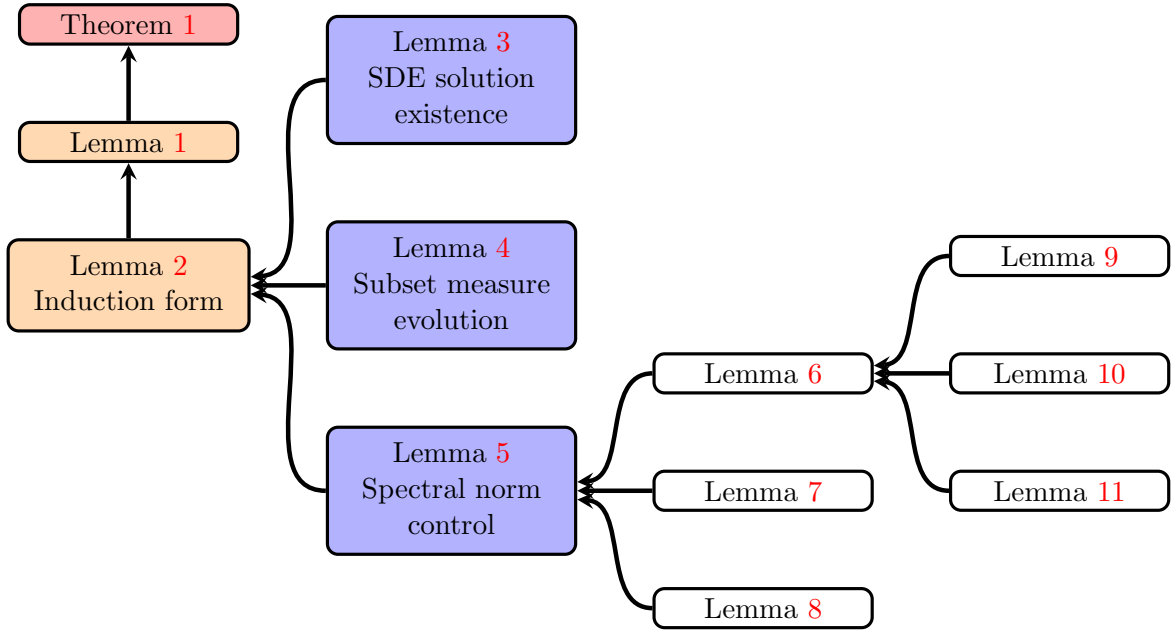
\begin{figure}[ht]
\begin{tikzpicture}[>=stealth,every node/.style={shape=rectangle,draw,rounded corners, minimum width=2.9cm,},]
    % create the nodes
    \node[very thick, fill=red!30] (t1) {Theorem~\ref{thm:main_KLS}};
    \node[very thick, fill=orange!30] (l1) [below=of
    t1]{Lemma~\ref{lem:main_lemma_KLS}} ;
    \node[very thick, fill=orange!30] (l2) [below=of l1]{\begin{tabular}{c}
      Lemma~\ref{lem:induction_isoperi} \\ Induction form
    \end{tabular}};
    \node[very thick, fill=blue!30, minimum width=3.3cm] (l4) [right=of l2]{
    \begin{tabular}{c}
      Lemma~\ref{lem:g_control} \\ Subset measure \\ evolution
    \end{tabular}};
    \node[very thick, fill=blue!30, minimum width=3.3cm] (l3) [above=of l4]{
    \begin{tabular}{c}
      Lemma~\ref{lem:existence_SDE_solution} \\ SDE solution \\ existence
    \end{tabular}};
    \node[very thick, fill=blue!30, minimum width=3.3cm] (l5) [below=of l4]{
    \begin{tabular}{c}
      Lemma~\ref{lem:isoperi_implies_A_control} \\ Spectral norm \\ control
    \end{tabular}};
    \node[very thick] (l7) [right=of l5]{Lemma~\ref{lem:poten_h_bound}};
    \node[very thick] (l6) [above=of l7]{Lemma~\ref{lem:potentwo_derivative_bound}};
    \node[very thick] (l8) [below=of l7]{Lemma~\ref{lem:poten_f_bound}};

    \node[very thick] (l10) [right=of l6]{Lemma~\ref{lem:ttensor_A_isoperimetric_bound}};
    \node[very thick] (l9) [above=of l10]{Lemma~\ref{lem:tensor_vector_bound}};
    \node[very thick] (l11) [below=of l10]{Lemma~\ref{lem:ttensor_oneoverT_bound}};

    % connect the nodes
    \draw[->, line width=.6mm] (l1) to[out=90,in=270](t1);
    \draw[->,  line width=.6mm] (l2) to[out=90,in=270] (l1);
    \draw[->,  line width=.6mm] (l3) to[out=180,in=4] (l2);
    \draw[->,  line width=.6mm] (l4) to[out=180,in=0] (l2);
    \draw[->,  line width=.6mm] (l5) to[out=180,in=-4] (l2);

    \draw[->,  line width=.6mm] (l6) to[out=180,in=4] (l5);
    \draw[->,  line width=.6mm] (l7) to[out=180,in=0] (l5);
    \draw[->,  line width=.6mm] (l8) to[out=180,in=-4] (l5);

    \draw[->,  line width=.6mm] (l9) to[out=180,in=4] (l6);
    \draw[->,  line width=.6mm] (l10) to[out=180,in=0] (l6);
    \draw[->,  line width=.6mm] (l11) to[out=180,in=-4] (l6);
    % \draw[->, dashed,  line width=0.6mm] (l1) to[out=0,in=180] (t2);
    % \draw[->, dashed,  line width=0.6mm] (l2) to[out=0,in=187] (t2);
    % \draw[->, dashed,  line width=0.6mm] (l4) to[out=180,in=195] (t2);
\end{tikzpicture}
  \caption{Proof sketch}
  \label{fig:proof_sketch}
\end{figure}

The same scheme was refined in Lee and Vempala~\cite{lee2016eldan} to obtain the current best lower bound in the KLS conjecture. Lee and Vempala directly attack the KLS conjecture while following the same stochastic localization scheme to smoothly multiply a Gaussian part to the log-concave density. Their use of a new potential function leads to the current best lower bound in the KLS conjecture. The proof in this paper builds on Lee and Vempala~\cite{lee2016eldan}'s refinements of Eldan's method, while it improves the handling of several quantities involved in the stochastic localization scheme. Figure~\ref{fig:proof_sketch} provides a diagram showing the relationship between the main lemmas.

To ensure the existence and the uniqueness of the stochastic localization construction, we first prove a lemma that deals with log-concave densities with compact support. Then we relate back to the main theorem by finding a compact support which contains most of the probability measure for a log-concave density.

\begin{lemma}
  \label{lem:main_lemma_KLS}
  There exists a universal constant $c$ such that for any log-concave density $p$ in $\real^\dims$ with compact support and any integer $\ell \geq 1$, we have
  \begin{align}
    \label{eq:main_lemma_KLS}
    \isoperi(p) \geq \frac{1}{\brackets{c \cdot \ell \parenth{\log(\dims)+1}}^{\ell/2} \dims^{16/\ell} \cdot \sqrt{\spec{p}}}.
  \end{align}
\end{lemma}

The proof of Lemma~\ref{lem:main_lemma_KLS} is provided in Section~\ref{sub:proof_of_lemma_lem:main_lemma_kls} after we introduce the intermediate lemmas. The use of the integer $l$ in the lemma indicates that we control the Cheeger isoperimetric coefficient in an iterative fashion. In fact, we prove Lemma~\ref{lem:main_lemma_KLS} by induction over $l$ starting from the known bound in Equation~\eqref{eq:KLS_original_bound}. For this, we define the supremum of the product of the isoperimetric coefficient and the square-root of its spectral norm over all log-concave densities in $\real^\dims$ with compact support:
\begin{align}
  \label{eq:def_inf_isoperi}
  \isoperi_\dims = \inf_{\text{ \shortstack{log-concave density $p$ in $\real^\dims$ \\ with compact support}}} \isoperi(p) \sqrt{\spec{p}}.
\end{align}
Then we prove the following lemma on the lower bound of $\isoperi_\dims$, which serves as the main induction argument.
\begin{lemma}
  \label{lem:induction_isoperi}
  Suppose that $\isoperi_k \geq \frac{1}{\alpha k^\beta}$ for all $k \leq \dims$ for some $0 \leq \beta \leq \frac{1}{2}$ and $\alpha \geq 1$, take $q = \lceil \frac{1}{\beta} \rceil + 1$, there exists a universal constant $c$ such that we have
  \begin{align*}
    \isoperi_\dims \geq \frac{1}{c \cdot q^{1/2} \alpha \log(\dims)^{1/2} \dims^{\beta - \beta / (8q) }}.
  \end{align*}
\end{lemma}
The proof of Lemma~\ref{lem:induction_isoperi} is provided towards the end of this section in Section~\ref{sub:proof_of_lemma_lem:induction_isoperi}. To have a good understanding of how we get there, we start by introducing the stochastic localization scheme introduced by Eldan~\cite{eldan2013thin}.

\subsection{Eldan's stochastic localization scheme} % (fold)
\label{sub:the_stochastic_localization_scheme}
Given a log-concave density $p$ in $\real^\dims$ with covariance matrix $\co$, we define the following stochastic differential equation (SDE)
\begin{align}
  \label{eq:def_SDE_ctBt}
  dc_t &= C_t^{1/2}dW_t + C_t \mu_t dt,\quad c_0 = 0,\\
  dB_t &= C_t dt,\quad B_0 = 0, \notag
\end{align}
where $W_t$ is the Wiener process, the matrix $C_t$, the density $p_t$, the mean $\mu_t$ and the covariance $A_t$ are defined as follows
\begin{align}
  \label{eq:def_SDE_parameters}
  C_t &= \co^{-1}, \\
  p_t(x) &= \frac{e^{c_t\tp x - \frac{1}{2}x\tp B_t x} p(x)}{\int_{\real^\dims} e^{c_t \tp x - \frac{1}{2}y\tp B_t y} p(y) dy}, \text{for $x \in \real^\dims$},\\
  \mu_t &= \int_{\real^\dims} x p_t(x)dx, \\
  A_t &= \int_{\real^\dims} \parenth{x - \mu_t}\parenth{x - \mu_t} p_t(x) dx,
\end{align}

The next lemma shows the existence and the uniqueness of the SDE solution.
\begin{lemma}
  \label{lem:existence_SDE_solution}
  Given a density $p$ in $\real^\dims$ with compact support with covariance $\co$ and $\co$ is invertible, then the SDE~\eqref{eq:def_SDE_ctBt} is well defined and it has a unique solution on the time interval $[0, T]$, for any time $T > 0$. Additionally, for any $x \in \real^\dims$, $p_t(x)$ is a martingale with
  \begin{align}
    \label{eq:derivative_pt}
    dp_t(x) = \parenth{x - \mu_t}\tp \co^{-1/2} dW_t p_t(x).
  \end{align}
\end{lemma}
The proof of Lemma~\ref{lem:existence_SDE_solution} follows from the standard existence and uniqueness theorem of SDE (Theorem 5.2 in {\O}ksendal~\cite{oksendal2003stochastic}). The proof is provided in Appendix~\ref{sec:proof_of_derivatives}.

Before we dive into the proof of Lemma~\ref{lem:induction_isoperi}, we discuss how the stochastic localization scheme allows us to control the boundary measure of a subset. First, according to the concavity of the isoperimetric profile (Theorem 2.8 in Sternberg and Zumbrun~\cite{sternberg1999connectivity} or Theorem 1.8 in Milman~\cite{milman2009role}), it is sufficient to consider subsets of measure $1/2$ in the definition of the isoperimetric coefficient in Equation~\eqref{eq:def_isoperi}. Second, the density $p_t$ is log-concave and it is more log-concave than the Gaussian density proportional to $e^{-\frac{1}{2}x\tp B_t x}$. It can be shown via the KLS localization lemma~\cite{kannan1995isoperimetric} that a density which is more log-concave than a Gaussian has an isoperimetric coefficient lower bound that depends on the covariance of the Gaussian (see e.g. Theorem 2.7 in Ledoux~\cite{ledoux2001concentration} or Theorem 4.4 in Cousins and Vempala~\cite{cousins2014cubic}).  Third, given an initial subset $E$ of $\real^\dims$ with measure $p(E) = \frac{1}{2}$, using the martingale property of $p_t(E)$, we observe that
\begin{align*}
  p(\partial E) &= \Exs \brackets{p_t(\partial E)} \\
  & \stackrel{(i)}{\geq} \Exs\brackets{\frac{1}{2}\vecnorm{B_t^{-1}}{2}^{-1/2} \min \parenth{p_t(E), p_t(E^c)}} \\
  & \stackrel{(ii)}{\geq} \frac{1}{4}\cdot \frac{1}{2}\vecnorm{B_t^{-1}}{2}^{-1/2} \Prob( \frac{1}{4}\leq p_t(E) \leq \frac{3}{4})\\
  & = \frac{1}{4}\vecnorm{B_t^{-1}}{2}^{-1/2} \Prob( \frac{1}{4}\leq p_t(E) \leq \frac{3}{4}) \cdot \min\braces{p(E), p(E^c)}.
\end{align*}
Inequality (i) uses the isoperimetric inequality for a log-concave density which is more log-concave than a Gaussian density proportional to $e^{-\frac{1}{2}x\tp B_t x}$ \cite{ledoux2001concentration,cousins2014cubic}. Inequality (ii) uses the fact that $p_t(E)$ is nonnegative.

Based on the above observation, the high level idea of the proof requires two main steps:
\begin{itemize}
  \item There exists some time $t > 0$, such that the Gaussian component $\frac{1}{2}x\tp B_t x$ of the density $p_t$ is large enough, so that we can apply the known isoperimetric inequality for densities more log-concave than a Gaussian.
  \item We need to control the quantity $p_t(E)$ so that the obtained isoperimetric inequality at time $t$ can be related back to that at time 0.
\end{itemize}
The first step is obvious since our construction explicitly enforces the density $p_t$ to have a Gaussian component $\frac{1}{2}x\tp B_t x$ in Equation~\eqref{eq:def_SDE_parameters}. Then the remaining question is whether we can run the SDE long enough to make the Gaussian component large enough while still keeping $p_t(E)$ to be the same order as $p(E) = \frac{1}{2}$ with large probability.

% subsection the_stochastic_localization_scheme (end)

\subsection{Control the evolution of the measure of a subset} % (fold)
\label{sub:control_the_evolution_of_the_measure_of_a_subset}

\begin{lemma}
  \label{lem:g_control}
  Under the same assumptions of Lemma~\ref{lem:existence_SDE_solution}, for any measurable subset $E$ of $\real^\dims$ with $p(E) = \frac{1}{2}$ and $t > 0$, the solution $p_t$ of the SDE~\eqref{eq:def_SDE_parameters} satisfies
  \begin{align*}
    \Prob\parenth{\frac{1}{4} \leq p_t(E) \leq \frac{3}{4}} \geq \frac{9}{10} - \Prob\parenth{\int_0^t \vecnorm{\co^{-1/2} A_t \co^{-1/2}}{2} ds \geq \frac{1}{64}}.
  \end{align*}
\end{lemma}
This lemma is proved in Lemma 29 of Lee and Vempala~\cite{lee2016eldan}. We provide a proof here for completeness.

\paragraph{Proof of Lemma~\ref{lem:g_control}:} Let $g_t = p_t(E)$.
Using Equation~\eqref{eq:derivative_pt}, we obtain the following derivative of $g_t$
\begin{align*}
  d g_t &= \int_E (x - \mu_t)\tp \co^{-1/2} dW_t p_t(x) dx.
\end{align*}
Its quadratic variation is
\begin{align*}
  d\brackets{g}_t &= \vecnorm{\int_E \co^{-1/2} (x - \mu_t) p_t(x) dx }{2}^2 dt \\
  &= \max_{\vecnorm{\xi}{2} \leq 1} \parenth{\int_E \xi\tp \co^{-1/2} (x - \mu_t) p_t(x) dx}^2 dt \\
  &\leq \max_{\vecnorm{\xi}{2} \leq 1} \parenth{\int_E \parenth{\xi\tp \co^{-1/2} (x - \mu_t)}^2 p_t(x) dx} \parenth{\int_E p_t(x) dx} dt \\
  &\leq \max_{\vecnorm{\xi}{2} \leq 1} \xi\tp \co^{-1/2} A_t \co^{-1/2} \xi dt \\
  &= \vecnorm{\co^{-1/2} A_t \co^{-1/2}}{2} dt,
\end{align*}
where the inequality follows from Cauchy-Schwarz inequality.
Applying the Dambis, Dubins-Schwarz theorem (see e.g. Revuz and Yor~\cite{revuz2013continuous} Section V.1 Theorem 1.7), there exists a Wiener process $\tilde{W}_t$ such that $g_t - g_0$ has the same distribution as $\tilde{W}_{[g]_t}$. Since $g_0 = \frac{1}{2}$, we obtain
\begin{align*}
  \Prob\parenth{\frac{1}{4} \leq p_t(E) \leq \frac{3}{4}} &= \Prob\parenth{-\frac{1}{4} \leq \tilde{W}_{[g]_t} \leq \frac{1}{4}} \\
  &\geq 1 - \Prob\parenth{\max_{0 \leq s \leq \frac{1}{64}} \abss{\tilde{W}_s} > \frac{1}{4}} - \Prob([g]_t > \frac{1}{64}) \\
  &= 1 - 4 \Prob\parenth{\tilde{W}_{\frac{1}{64}} > \frac{1}{4}} - \Prob\parenth{[g]_t > \frac{1}{64}} \\
  &\geq \frac{9}{10} - \Prob\parenth{\int_0^t \vecnorm{\co^{-1/2} A_t \co^{-1/2}}{2} ds > \frac{1}{64}},
\end{align*}
where the last inequality follows from the fact that $\Prob\parenth{\xi > 2} < 0.023$ for $\xi$ follows the standard Gaussian distribution.

% \begin{lemma}
%   \label{lem:A_control_implies_isoperi}
%   For a logconcave density $p$ in $\real^\dims$ with compact support. If there exists $T>0$ such that the solution density $p_t$ under the stochastic localization scheme with initialization $p$ satisfies
%   \begin{align*}
%     \Prob\parenth{\int_{0}^T \vecnorm{C_t^{1/2}A_t C_t^{1/2}}{2} dt \geq c} \leq c,
%   \end{align*}
%   then we have
%   \begin{align*}
%     \isoperi(p) \leq T^{-1/2}
%   \end{align*}
% \end{lemma}
% \paragraph{Proof of Lemma~\ref{lem:A_control_implies_isoperi}:}
% It is sufficient to consider all subsets of measure $\frac{1}{2}$.  Consider any measurable subset $E$ of $\real^\dims$ with $p(E) = \frac{1}{2}$. Since $p_t$ is a martingale, we have for any $t > 0$,
% \begin{align*}
%   p(\partial E) = \Exs \brackets{ p_t \parenth{\partial E}}.
% \end{align*}
% Since $p_{T}$ has a Gaussian component, we can apply the isoperimetic inequality for densities with a Gaussian component and transfer this property to the initial density $p$.
% \begin{align*}
%   p(\partial E) &= \Exs \brackets{ p_{T} \parenth{\partial E}} \\
%   &\stackrel{(i)}{=} \frac{1}{T^{1/2}} \Exs \brackets{\min\braces{g_T(E), 1- g_T(E^c)}} \\
%   &\geq \frac{1}{4 T^{1/2}} \Prob \parenth{ \frac{1}{4} \leq g_T(E) \leq \frac{3}{4}} \\
%   &\geq \frac{1}{4 T^{1/2}} \brackets{0.9 - \Prob\parenth{\int_{0}^T \vecnorm{C_t^{1/2}A_t C_t^{1/2}}{2} dt \geq \frac{1}{64}}}
% \end{align*}

% It remains to control the spectral norm of $C_t^{-1/2} A_t C_t^{-1/2}$.

% subsection controlling_the_mass_of_a_specific_set_e_at_time_t (end)

\subsection{Control the evolution of the spectral norm} % (fold)
\label{sub:control_the_evolution_of_the_spectral norm}
According to Lemma~\ref{lem:g_control}, to control the evolution of the measures of subsets, we need to control the spectral norm of $\co^{-1/2} A_t \co^{-1/2}$. The following lemma serves the purpose.
\begin{lemma}
  \label{lem:isoperi_implies_A_control}
  In addition to the same assumptions of Lemma~\ref{lem:existence_SDE_solution}, if $\isoperi_k \geq \frac{1}{\alpha k^\beta}$ for all $k \leq \dims$ for some $0 < \beta \leq \frac{1}{2}$ and $\alpha \geq 1$, then there exists a universal constant $c$ such that for $q = \lceil \frac{1}{\beta} \rceil + 1$, $\dims \geq 3$ and $T_2 = \frac{1}{ c \cdot q \alpha^2\log(\dims) \dims^{2\beta - \beta/(4q)}}$, we have
  \begin{align*}
    \Prob\parenth{\int_{0}^{T_2} \vecnorm{\co^{-1/2} A_t \co^{-1/2}}{2} dt \geq \frac{1}{64}} < \frac{4}{10}.
  \end{align*}
\end{lemma}
Direct control of the largest eigenvalue of $\co^{-1/2} A_t \co^{-1/2}$ is not trivial, instead we use the potential function $\potentwo_t$ to upper bound the largest eigenvalue.
Define
\begin{align}
  \label{eq:poten2_def}
  \qo_t &= \co^{-1/2} A_t \co^{-1/2} \notag \\
  \potentwo_t &= \trace\parenth{\qo_t^q}.
\end{align}
It is clear that $\potentwo_t^{1/q} \geq \vecnorm{\co^{-1/2} A_t \co^{-1/2}}{2}$. So in order to upper bound $\vecnorm{\co^{-1/2} A_t \co^{-1/2}}{2}$, it is sufficient to upper bound $\potentwo_t^{1/q}$. The advantage of using $\potentwo_t$ is that it is differentiable. We have the following differential for $A_t$ and $\potentwo_t$:
\begin{align}
  \label{eq:known_SDE_change_At}
  dA_t &= \int (x - \mu_t) (x - \mu_t)\tp \parenth{(x-\mu_t)\tp \co^{-1/2}dW_t} p_t(x) dx - A_t\co^{-1}A_t dt, \\
  \label{eq:poten2_derivative}
  d\potentwo_t &=  q \int \parenth{x-\mu_t}\tp \co^{-1/2} \parenth{\qo_t}^{q-1} \co^{-1/2} \parenth{x-\mu_t} \parenth{x-\mu_t}\tp \co^{-1/2} dW_t p_t(x) dx \notag \\
  &- q \trace\parenth{\qo_t^{q+1} } dt \notag \\
  &+ \frac{q}{2} \sum_{a = 0}^{q-2} \int \int \parenth{x-\mu_t}\tp \co^{-1/2} \qo_t^{a} \co^{-1/2} \parenth{y-\mu_t} \notag \\
  &\cdot \parenth{x-\mu_t}\tp \co^{-1/2} \qo_t^{q-2-a} \co^{-1/2} \parenth{y-\mu_t} \parenth{x-\mu_t}\tp \co^{-1} \parenth{y-\mu_t} p_t(x) p_t(y) dx dy dt.
\end{align}
Obtaining these differentials uses It\^o's formula and the proofs are provided in Appendix~\ref{sec:proof_of_derivatives}.

The next lemma upper bounds the terms in the potential $\potentwo_t$.

\begin{lemma}
  \label{lem:potentwo_derivative_bound}
  Under the same assumptions of Lemma~\ref{lem:isoperi_implies_A_control}, the potential $\potentwo_t$ defined in Equation~\eqref{eq:poten2_def} can be written as follows
  \begin{align*}
    d\potentwo_t = v_t\tp dW_t + \delta_t dt,
  \end{align*}
  where $v_t \in \real^\dims$ and $\delta_t \in \real$ satisfy
  \begin{align*}
    \vecnorm{v_t}{2} &\leq 16 q \potentwo_t^{1 + 1/(2q)}, \text{ and } \\
    \delta_t &\leq \min\braces{64 q^2 \alpha^2 \log(\dims) \dims^{2\beta -1/q}\potentwo_t^{1 + 1/q}, \frac{2q^2}{t} \potentwo_t}.
  \end{align*}
\end{lemma}
The proof of Lemma~\ref{lem:potentwo_derivative_bound} is provided in Section~\ref{sub:tensor_bounds}. Remark that bounds similar to the first bound of $\delta_t$ in Lemma~\ref{lem:potentwo_derivative_bound} have appeared in Lee and Vempala~\cite{lee2016eldan}, whereas the second bound of $\delta_t$ in Lemma~\ref{lem:potentwo_derivative_bound} is novel. The second bound of $\delta_t$ also leads to the following Lemma~\ref{lem:poten_f_bound} which gives better control of the potential than the previous proof by Lee and Vempala~\cite{lee2016eldan} when $t$ is large.

Using the bounds in Lemma~\ref{lem:potentwo_derivative_bound}, we state the two lemmas which control the potential $\potentwo_t$ in two ways.
\begin{lemma}
  \label{lem:poten_h_bound}
  Under the same assumptions of Lemma~\ref{lem:potentwo_derivative_bound}, using the following transformation
  \begin{align*}
    h: \real_+ &\rightarrow \real \\
    a &\mapsto -(a+1)^{-1/q}
  \end{align*}
  we have
  \begin{align*}
    \Prob\parenth{\max_{t \in [0, T_1]} h(\potentwo_t) \geq - \frac{1}{2}\parenth{\dims+1}^{-1/q} } \leq \exp(-\frac{2}{3} q\log(\dims)) \leq \frac{3}{10}
  \end{align*}
  where $T_1 = \frac{1}{32768 q \alpha^2 \log(\dims) \dims^{2\beta}}$.
\end{lemma}

\begin{lemma}
  \label{lem:poten_f_bound}
  Under the same assumptions of Lemma~\ref{lem:potentwo_derivative_bound}, using the following transformation
  \begin{align*}
    f: \real_+ &\rightarrow \real \\
    a &\mapsto a^{1/q}
  \end{align*}
  we have
  \begin{align*}
    \Exs f(\potentwo_{t_2}) \leq \Exs f(\potentwo_{t_1}) \parenth{\frac{t_2}{t_1}}^{2q}, \forall t_2 > t_1 > 0.
  \end{align*}
\end{lemma}
The proofs of Lemma~\ref{lem:poten_h_bound} and~\ref{lem:poten_f_bound} are provided in Section~\ref{sub:control_the_potential_in_two_time_periods}.

Now we are ready to prove Lemma~\ref{lem:isoperi_implies_A_control}.
\paragraph{Proof of Lemma~\ref{lem:isoperi_implies_A_control}:} We take
\begin{align*}
  T_1 = \frac{1}{32768 q \alpha^2 \log(\dims) \dims^{2\beta}}, \quad T_2 = \frac{\dims^{\beta/(4q)}}{40} T_1 = \frac{1}{ 1310720 q \alpha^2\log(\dims) \dims^{2\beta - \beta/(4q)}}.
\end{align*}
We bound the spectral norm of $\co^{-1/2}A_t \co^{-1/2}$ in two time intervals via Lemma~\ref{lem:poten_h_bound} and Lemma~\ref{lem:poten_f_bound}.
In the first time interval $[0, T_1]$, we have
\begin{align}
  \label{eq:spec_norm_first_time_period}
  \Prob\parenth{\int_0^{T_1} \vecnorm{\co^{-1/2}A_t \co^{-1/2}}{2} \geq \frac{1}{128}} &\leq \Prob\parenth{\max_{t \in [0, T_1]} \vecnorm{\co^{-1/2}A_t \co^{-1/2}}{2} \geq \frac{1}{128T_1}} \notag \\
  &\stackrel{(i)}{\leq} \Prob\parenth{\max_{t \in [0, T_1]} \vecnorm{\co^{-1/2}A_t \co^{-1/2}}{2} \geq 3 \dims^{1/q}} \notag \\
  &\stackrel{(ii)}{\leq} \Prob\parenth{\max_{t \in [0, T_1]} \potentwo_t \geq 3^{q} \dims} \notag \\
  &\stackrel{(iii)}{\leq} \Prob\parenth{\max_{t \in [0, T_1]} \potentwo_t + 1 \geq 2^{q} (\dims+1)} \notag \\
  &\stackrel{\phantom{(iii)}}{=} \Prob\parenth{\max_{t \in [0, T_1]} h(\potentwo_t) \geq -\frac{1}{2} \parenth{\dims+1}^{-1/q}} \notag \\
  &\stackrel{(iv)}{\leq} \frac{3}{10}.
\end{align}
Inequality (i) follows from the condition $\beta q \geq 1$. (ii) follows from the fact that $\trace\parenth{A^q}^{1/q} \geq \vecnorm{A}{2}$. (iii) is because $3^q \dims \geq 2^q (\dims + 1)$ when $q \geq 2$ and $\dims \geq 1$. $h$ is defined in Lemma~\ref{lem:poten_h_bound}. (iv) follows from Lemma~\ref{lem:poten_h_bound}.

In the first time interval, we can also bound the expectation of $\potentwo_{T_1}^{1/q}$. Since the density $p_{T_1}$ is more log-concave than a Gaussian density with covariance matrix $\frac{\co}{T_1}$, the covariance matrix of $p_{T_1}$ is upper bounded as follows (see Theorem 4.1 in Brascamp-Lieb~\cite{brascamp2002extensions} or Lemma 5 in Eldan and Lehec~\cite{eldan2014bounding})
\begin{align}
  \label{eq:covAt_bound}
  A_{T_1} \preceq \frac{\co}{T_1}.
\end{align}
Consequently, all the eigenvalues of $Q_{T_1}$ are less than $\frac{1}{T_1}$ and $\potentwo_{T_1}$ is upper bounded by $\frac{\dims}{T_1^{q}}$. Using the above bound, we can bound the expectation of $\potentwo_{T_1}^{1/q}$ as follows
\begin{align}
  \label{eq:exp_spec_norm_T1_bound}
  \Exs\brackets{\potentwo_{T_1}^{1/q}} &= \Exs\brackets{\indicator_{\potentwo_{T_1} \geq 3^q \dims} \potentwo_{T_1}^{1/q} + \indicator_{\potentwo_{T_1} < 3^q \dims} \potentwo_{T_1}^{1/q} } \notag \\
  &\stackrel{(i)}{\leq} \frac{\dims^{1/q}}{T_1} \exp\parenth{- \frac{2}{3} q \log(\dims)} + 3 \dims^{1/q}  \notag \\
  &\stackrel{(ii)}{\leq} 32768 \dims^{1/q} q \alpha^2 + 4 \dims^{1/q}  \notag \\
  &\leq 40000 \dims^{1/q} q \alpha^2.
\end{align}
Inequality (i) follows from Lemma~\ref{lem:poten_h_bound}, the inequality $3^q\dims \geq 2^q(\dims+1)$ (similar to what we did in the last four steps of Equation~\eqref{eq:spec_norm_first_time_period}) and Equation~\eqref{eq:covAt_bound}. (ii) follows from $q \geq 2$, $\beta \leq {1/2}$ and $\dims^{1/2} \geq \log(\dims)$ for $\dims \geq 3$.

In the second time interval, for $t \in [T_1, T_2]$, we have
\begin{align}
  \label{eq:exp_spec_norm_T2_bound}
  \Exs\brackets{\vecnorm{\co^{-1/2} A_{t} \co^{-1/2}}{2}} &\leq \Exs\brackets{\potentwo_{t}^{1/q}} \notag \\
  &\stackrel{(i)}{\leq} \Exs\brackets{\potentwo_{T_1}^{1/q}} \parenth{\frac{t}{T_1}}^{2q} \notag \\
  &\stackrel{(ii)}{\leq} \Exs\brackets{\potentwo_{T_1}^{1/q}} \parenth{\frac{T_2}{T_1}}^{2q} \notag \\
  &\stackrel{(iii)}{\leq} 1000 \dims^{\beta/2 + 1/q} q \alpha^2
\end{align}
Inequality (i) follows from Lemma~\ref{lem:poten_f_bound}. (ii) is because $t \leq T_2$. (iii) follows from $T_2 = \frac{\dims^{\beta/(4q)}}{40} T_1$.
Using the above bound, we control the spectral norm in the second time interval via Markov's inequality
\begin{align}
  \label{eq:spec_norm_second_time_period}
  \Prob\parenth{\int_{T_1}^{T_2} \vecnorm{\co^{-1/2}A_t \co^{-1/2}}{2} \geq \frac{1}{128}}
  &\stackrel{(i)}{\leq} \frac{\Exs \brackets{\int_{T_1}^{T_2} \vecnorm{\co^{-1/2}A_t \co^{-1/2}}{2} dt } }{1/128} \notag \\
  &\stackrel{(ii)}{\leq} T_2 \cdot 1000 \dims^{\beta/2 + 1/q} q \alpha^2 \cdot 128 \notag \\
  &\stackrel{(iii)}{<} \frac{1}{10},
\end{align}
where inequality (i) follows from Markov's inequality and (ii) follows from Equation~\eqref{eq:exp_spec_norm_T2_bound}. (iii) follows from the definition of $T_2$ and $\frac{\beta}{2}+\frac{1}{q} \leq 2\beta-\beta/(4q)$ when $\beta q \geq 1$ and $q \geq 2$.

Combining the bounds in the first and second time intervals in Equation~\eqref{eq:spec_norm_first_time_period} and~\eqref{eq:spec_norm_second_time_period}, we obtain
\begin{align}
  &\Prob\parenth{\int_{0}^{T_2} \vecnorm{\co^{-1/2}A_t \co^{-1/2}}{2} \geq \frac{1}{64}} \notag \\
  \leq & \Prob\parenth{\int_{0}^{T_1} \vecnorm{\co^{-1/2}A_t \co^{-1/2}}{2} \geq \frac{1}{128}} + \Prob\parenth{\int_{T_1}^{T_2} \vecnorm{\co^{-1/2}A_t \co^{-1/2}}{2} \geq \frac{1}{128}} \leq \frac{4}{10}.
\end{align}

% subsection controlling_the_evolution_of_the_largest_eigenvalue (end)

\subsection{Proof of Lemma~\ref{lem:induction_isoperi}} % (fold)
\label{sub:proof_of_lemma_lem:induction_isoperi}
The proof of Lemma~\ref{lem:induction_isoperi} follows the strategy described after Lemma~\ref{lem:existence_SDE_solution}. We make the arguments rigorous here. We consider a log-concave density $p$ in $\real^\dims$ with compact support. Without loss of generality, we can assume that the covariance matrix $A$ of the density $p$ is invertible. Otherwise, the density $p$ is degenerate and we can instead prove the results in a lower dimension.

According to the concavity of the isoperimetric profile (Theorem 2.8 in Sternberg and Zumbrun~\cite{sternberg1999connectivity} or Theorem 1.8 in Milman~\cite{milman2009role}), it is sufficient to consider subsets of measure $1/2$ in the definition of isoperimetric coefficient~\eqref{eq:def_isoperi}. Given an initial subset $E$ of $\real^\dims$ with $p(E) = \frac{1}{2}$, use the martingale property of $p_{T_2}(E)$, we have
\begin{align*}
  p(\partial E) &= \Exs \brackets{p_{T_2}(\partial E)} \\
  & \stackrel{(i)}{\geq} \Exs\brackets{\frac{1}{2}\vecnorm{B_{T_2}^{-1}}{2}^{-1/2} \min \parenth{p_{T_2}(E), p_{T_2}(E^c)}} \\
  & \stackrel{(ii)}{\geq} \frac{1}{4}\cdot \frac{1}{2}\vecnorm{B_{T_2}^{-1}}{2}^{-1/2} \Prob( \frac{1}{4}\leq p_{T_2}(E) \leq \frac{3}{4})\\
  & = \frac{1}{4}\vecnorm{B_{T_2}^{-1}}{2}^{-1/2} \Prob( \frac{1}{4}\leq p_{T_2}(E) \leq \frac{3}{4}) \cdot \min\braces{p(E), p(E^c)} \\
  &\stackrel{(iii)}{\geq} \frac{1}{8}\vecnorm{B_{T_2}^{-1}}{2}^{-1/2} \cdot \min\braces{p(E), p(E^c)} \\
  &\stackrel{(iv)}{=} \frac{1}{8}T_2^{1/2}\vecnorm{\co}{2}^{-1/2} \cdot \min\braces{p(E), p(E^c)}.
\end{align*}
Inequality (i) uses the isoperimetric inequality for a log-concave density which is more log-concave than a Gaussian density proportional to $e^{-\frac{1}{2}x\tp B_t x}$ (see e.g. Theorem 2.7 in Ledoux~\cite{ledoux2001concentration} or Theorem 4.4 in Cousins and Vempala~\cite{cousins2014cubic}). Inequality (ii) follows from the fact that $p_t(E)$ is nonnegative. (iii) follows from Lemma~\ref{lem:g_control} and Lemma~\ref{lem:isoperi_implies_A_control} (for $\dims \geq 3$). (iv) follows from the construction that $B_t = t \co^{-1}$.
We conclude the proof since $T_2$ is taken as $\frac{1}{ c \cdot q \alpha^2\log(\dims) \dims^{2\beta - \beta/(4q)}}$ with $c$ as a constant.
The above proof only works for $\dims \geq 3$. It is easy to verify that Lemma~\ref{lem:induction_isoperi} still holds for the case for $\dims = 1, 2$ from the original KLS bound in Equation~\eqref{eq:KLS_original_bound}.

% subsection proof_of_lemma_lem:induction_isoperi (end)

\subsection{Proof of Lemma~\ref{lem:main_lemma_KLS}} % (fold)
\label{sub:proof_of_lemma_lem:main_lemma_kls}
The proof of Lemma~\ref{lem:main_lemma_KLS} consists of applying Lemma~\ref{lem:induction_isoperi} recursively. We define
\begin{align*}
  \alpha_1 = 4, \beta_1 = \frac{1}{2}.
\end{align*}
For $\ell \geq 1$, we define $\alpha_\ell$ and $\beta_\ell$ recursively as follows:
\begin{align}
  \label{eq:alpha_beta_l_recursion}
  \alpha_{\ell+1} &= 2c \cdot \alpha_\ell \beta_\ell^{-1/2},  \notag \\
  \beta_{\ell+1} &= \beta_\ell - \beta_\ell^2/16,
\end{align}
where $c$ is the constant in Lemma~\ref{lem:induction_isoperi}.
It is not difficult to show by induction that $\alpha_\ell$ and $\beta_\ell$ satisfy
\begin{align}
  \label{eq:alpha_beta_l_bound}
  \frac{1}{\ell+1} &\leq \beta_\ell \leq \frac{16}{\ell} \notag \\
  \alpha_\ell &\leq \parenth{4c^2 \ell}^{\ell/2}.
\end{align}

We start with a known bound from the original KLS paper~\cite{kannan1995isoperimetric}
\begin{align*}
  \isoperi_\dims \geq \frac{1}{\alpha_1 \dims^{\beta_1}},\quad \forall \dims \geq 1.
\end{align*}
In the induction, suppose that we have
\begin{align*}
  \isoperi_\dims \geq \frac{1}{\alpha_\ell \parenth{\log(\dims)+1}^{\ell/2} \dims^{\beta_\ell}},\quad \forall \dims \geq 1.
\end{align*}
From the above inequality, we obtain for any $1 \leq k \leq \dims$,
\begin{align*}
  \isoperi_k \geq \frac{1}{\alpha_\ell' k^\beta_\ell},
\end{align*}
with $\alpha_\ell' = \alpha_\ell \parenth{\log(\dims) +1}^{\ell/2}$. Using the above lower bounds for $\isoperi_k$, we can apply Lemma~\ref{lem:induction_isoperi}. For integer $\ell+1$, we have
\begin{align*}
  \isoperi_\dims &\stackrel{(i)}{\geq} \frac{1}{c \cdot q^{1/2} \alpha_\ell \parenth{\log(\dims)+1}^{l/2} \log(\dims)^{1/2} \dims^{\beta_\ell - \beta_\ell / (8q) }}\\
  &\stackrel{(ii)}{\geq} \frac{1}{2c \cdot \alpha_\ell \beta_\ell^{-1/2}  \parenth{\log(\dims)+1}^{(l+1)/2} \dims^{\beta_\ell - \beta_\ell^2 / 16 }} \\
  &= \frac{1}{\alpha_{\ell+1} \parenth{\log(\dims)+1}^{(\ell+1)/2} \dims^{\beta_{\ell+1}}}
\end{align*}
where inequality (i) follows from Lemma~\ref{lem:induction_isoperi}, inequality (ii) follows from $q \leq \frac{2}{\beta}$ and the last equality follows from the definition of $\alpha_\ell$ and $\beta_\ell$. We conclude Lemma~\ref{lem:main_lemma_KLS} using the $\alpha_\ell$ and $\beta_\ell$ bounds in Equation~\eqref{eq:alpha_beta_l_bound}.

% subsection proof_of_lemma_lem:main_lemma_kls (end)

\subsection{Proof of Theorem~\ref{thm:main_KLS}} % (fold)
\label{sub:proof_of_theorem_thm:main_kls}
To derive Theorem~\ref{thm:main_KLS} from Lemma~\ref{lem:main_lemma_KLS}, it is sufficient to show that for any log-concave density $p$ in $\real^\dims$, most of its probability measure is on a compact support. Let $\mu$ be the mean of the density $p$. Since $r \mapsto p(\mathbb{B}\parenth{\mu, r}^c)$ is an non-increasing function of $r$ with limit $0$ at $\infty$, there exists a radius $R > 0$, such that $p(\mathbb{B}\parenth{\mu, R}^c) \leq 0.2$. Note that it is possible to get a better bound via e.g. log-concave concentration bounds from Paouris~\cite{paouris2006concentration}, but knowing the existence of such radius $R$ is sufficient for the proof here.

Denote $B = \mathbb{B}\parenth{\mu, R}$. Then $p(B^c)\leq 0.2$. Let $\varrho$ be the density obtained by truncating $p$ on the ball $B$. Then $\varrho$ is log-concave and it has compact support.
For a subset $E \subset \real^\dims$ of measure such that $p(E) = \frac{1}{2}$, we have
\begin{align*}
  p(\partial E) &\geq \varrho(\partial E) p(B) \\
  &\geq \isoperi(\varrho) \min\parenth{\varrho(E), \varrho(E^c)} p(B) \\
  &= \isoperi(\varrho) \min\parenth{p(E \cap B), p(B \cap E^c)} \\
  &\geq \isoperi(\varrho) \min\parenth{p(E) - p(B^c), p(E^c) - p(B^c)} \\
  &\geq \frac{1}{2} \isoperi(\varrho) \min\parenth{p(E), p(E^c)}.
\end{align*}
The last inequality follows because $p(E^c) - p(B^c) \geq 0.5 - 0.2 \geq \frac{1}{4}$.
Since it is sufficient to consider subsets of measure 1/2 in the definition of the isoperimetric coefficient~\cite{sternberg1999connectivity,milman2009role}, we conclude that the isoperimetric coefficient of $p$ is lower bounded by half of that of $\varrho$. Applying Lemma~\ref{lem:main_lemma_KLS} for the isoperimetric coefficient of $\varrho$, we obtain Theorem~\ref{thm:main_KLS}.

% subsection proof_of_theorem_thm:main_kls (end)

\section{Proof of auxiliary lemmas} % (fold)
\label{sec:proof_of_auxiliary_lemmas}

In this section, we prove auxiliary lemmas~\ref{lem:potentwo_derivative_bound},~\ref{lem:poten_h_bound} and~\ref{lem:poten_f_bound}.

\subsection{Tensor bounds and proof of Lemma~\ref{lem:potentwo_derivative_bound}} % (fold)
\label{sub:tensor_bounds}
In this subsection, we prove Lemma~\ref{lem:potentwo_derivative_bound}. Since Lemma~\ref{lem:potentwo_derivative_bound} involves the third-order moment tensor of a log-concave density, we define the following 3-Tensor for any probability density $p \in \real^\dims$ with mean $\mu$ to simplify notations.
\begin{align}
  \label{def:ttensor}
  \ttensor_p: \quad
  &\real^{\dims \times \dims} \times \real^{\dims \times \dims} \times \real^{\dims \times \dims} \rightarrow \real \notag \\
  &(A, B, C) \mapsto \int \int (x-\mu)\tp A (y-\mu) \cdot (x-\mu)\tp B (y - \mu) \cdot (x - \mu) \tp C (y - \mu) p(x) p(y) dx dy.
\end{align}
For $A, B, C$ three matrices in $\real^{\dims \times \dims}$, we can write $\ttensor_p(A, B, C)$ equivalently as
\begin{align*}
  \ttensor_p(A, B, C) = \Exs_{X, Y \sim p} (X-\mu) \tp A (Y-\mu) \cdot (X-\mu) \tp B (Y-\mu) \cdot (X-\mu) \tp C (Y-\mu).
\end{align*}
Before we prove Lemma~\ref{lem:potentwo_derivative_bound}, we prove the following properties related to the 3-Tensor.
\begin{lemma}
  \label{lem:tensor_vector_bound}
  Suppose $p$ is a log-concave density with mean $\mu$ and covariance $A$. Then for any positive semi-definite matrices $B$ and $C$, we have
  \begin{align*}
    \vecnorm{\int B^{1/2} (x - \mu) (x - \mu) \tp C (x - \mu) p(x)dx}{2} \leq 16 \vecnorm{A^{1/2}B A^{1/2}}{2}^{1/2} \trace\parenth{A^{1/2} C A^{1/2}}
  \end{align*}
\end{lemma}

\begin{lemma}
  \label{lem:ttensor_A_isoperimetric_bound}
  Suppose that $\isoperi_k \geq \frac{1}{\alpha k^\beta}$ for all $k \leq \dims$ for some $0 \leq \beta \leq \frac{1}{2}$ and $\alpha \geq 1$. Suppose $p$ is a log-concave density in $\real^\dims$ with covariance $A$ and $A$ is invertible. Then for $q \geq \frac{1}{2\beta}$, we have
  \begin{align*}
    \ttensor_p(A^{q-2}, \Ind_\dims, \Ind_\dims) \leq 128 \alpha^2 \log(\dims) \dims^{2\beta - 1/q} \trace(A^q) ^{1 + 1/q}.
  \end{align*}
\end{lemma}

\begin{lemma}
  \label{lem:ttensor_oneoverT_bound}
  Given $\tau > 0$. Suppose $p$ is a log-concave density which is more log-concave than $\Normal(0, \frac{1}{\tau} \Ind_\dims)$. Let $A$ be its covariance matrix. Suppose $A$ is invertible then for $q \geq 3$, we have
  \begin{align*}
    \ttensor_p(A^{q-2}, \Ind_\dims, \Ind_\dims) \leq \frac{4}{\tau} \trace\parenth{A^{q}}
  \end{align*}
\end{lemma}

\begin{lemma}
  \label{lem:ttensor_swap}
  Suppose $p$ is a log-concave density in $\real^\dims$. For any $\delta \in [0, 1]$, for $A, B, C$ positive semi-definite matrices then
  \begin{align}
    \label{eq:ttensor_swap}
    \ttensor_{p}(B^{1/2}A^\delta B^{1/2}, B^{1/2}A^{1-\delta}B^{1/2}, C) \leq \ttensor_{p}(B^{1/2}AB^{1/2}, B, C).
  \end{align}
\end{lemma}
The proofs of the above lemmas are provided in Section~\ref{sub:proof_of_tensor_bounds}.

Now we are ready to prove Lemma~\ref{lem:potentwo_derivative_bound}.
\paragraph{Proof of Lemma~\ref{lem:potentwo_derivative_bound}:} We first prove the bound on $\vecnorm{v_t}{2}$, where
\begin{align*}
  v_t = q \int \co^{-1/2} \parenth{x-\mu_t} \parenth{x-\mu_t}\tp \co^{-1/2} \parenth{\qo_t}^{q-1} \co^{-1/2} \parenth{x-\mu_t} p_t(x) dx.
\end{align*}
Applying Lemma~\ref{lem:tensor_vector_bound} and knowing the covariance of $p_t$ is $A_t$, we obtain
\begin{align*}
  \vecnorm{v_t}{2} &\leq 16 q \vecnorm{A_t^{1/2} \co^{-1} A_t^{1/2}}{2}^{1/2} \trace\parenth{A_t^{1/2} \co^{-1/2} Q_t^{q-1} \co^{-1/2} A_t^{1/2}} \\
  &\stackrel{(i)}{=} 16 q \vecnorm{A_t^{1/2} \co^{-1} A_t^{1/2}}{2}^{1/2} \trace\parenth{Q_t^{q}} \\
  &\stackrel{(ii)}{=} 16 q \vecnorm{Q_t}{2}^{1/2} \trace\parenth{Q_t^{q}} \\
  &\stackrel{(iii)}{\leq}  16 q \brackets{\trace\parenth{Q_t^{q}}}^{1+1/(2q)}.
\end{align*}
Equality (i) uses the definition of $Q_t = \co^{-1/2} A_t \co^{-1/2}$. Equality (ii) uses the fact that $\vecnorm{MM\tp}{2} = \vecnorm{M\tp M}{2}$ for any square matrix $M \in \real^{\dims \times \dims}$.
Inequality (iii) uses that $\vecnorm{M}{2} \leq \trace\parenth{M^q}^{1/q}$ for any positive semi-definite matrix $M$.

Next, we bound $\delta_t$ in two ways. We can ignore the negative term in $\delta_t$ to obtain the following:
\begin{align}
  \label{eq:delta_t_bound_intermediate}
  \delta_t &\leq \frac{q}{2} \sum_{a = 0}^{q-2} \int \int \parenth{x-\mu_t}\tp \co^{-1/2} \qo_t^{a} \co^{-1/2} \parenth{y-\mu_t} \notag \\
  &\cdot \parenth{x-\mu_t}\tp \co^{-1/2} \qo_t^{q-2-a} \co^{-1/2} \parenth{y-\mu_t} \parenth{x-\mu_t}\tp \co^{-1} \parenth{y-\mu_t} p_t(x) p_t(y) dx dy \notag \\
  &= \frac{q}{2} \sum_{a = 0}^{q-2} \ttensor_{\denb_t}(Q_t^{a}, Q_t^{q-2-a}, \Ind_\dims),
\end{align}
where $\denb_t$ is the density of linear-transformed random variable $\co^{-1/2}\parenth{X-\mu_t}$ for $X$ drawn from $p_t$ and $\mu_t$ is the mean of $p_t$. $\denb_t$ is still log-concave since any linear transformation of a log-concave density is log-concave (see e.g. Saumard and Wellner~\cite{saumard2014log}). $\denb_t$ has covariance $\co^{-1/2} A_t \co^{-1/2}$, which is also $Q_t$. For $a \in \braces{0, \cdots, q-2}$, we have
\begin{align*}
  \ttensor_{\denb_t}(\qo_t^{a}, \qo_t^{q-2-a}, \Ind_\dims) &\stackrel{(i)}{\leq} \ttensor_{\denb_t}(\qo_t^{q-2}, \Ind_\dims, \Ind_\dims) \\
  &\stackrel{(ii)}{\leq} 128 \alpha^2 \log(\dims) \dims^{2\beta - 1/q} \brackets{\trace\parenth{\qo_t^q}}^{1+1/q}.
\end{align*}
Inequality (i) follows from Lemma~\ref{lem:ttensor_swap}. Inequality (ii) follows from Lemma~\ref{lem:ttensor_A_isoperimetric_bound}. Since there are $q-1$ terms in the sum, we conclude the first part of the bound for $\delta_t$.

On the other hand, since $p_t$ is more log-concave than the Gaussian density proportional to $e^{-\frac{t}{2} (x-\mu_t)\tp \co^{-1} (x-\mu_t)}$, $\denb_t$ is more log-concave than the Gaussian density proportional to $e^{-\frac{t}{2} x\tp x}$. Applying Lemma~\ref{lem:ttensor_swap} and Lemma~\ref{lem:ttensor_oneoverT_bound} to each term in Equation~\eqref{eq:delta_t_bound_intermediate}, we obtain
\begin{align*}
  \delta_t &\leq \frac{q^2}{2} \ttensor_{\denb_t}(\qo_t^{q-2}, \Ind_\dims, \Ind_\dims) \\
  &\leq \frac{2q^2}{t} \trace\parenth{\qo_t^{q}}.
\end{align*}
This concludes the second part of the bound for $\delta_t$.
% subsection tensor_bounds (end)

\subsection{Control of the potential in two time intervals} % (fold)
\label{sub:control_the_potential_in_two_time_periods}
In this subsection, we prove Lemma~\ref{lem:poten_h_bound} and Lemma~\ref{lem:poten_f_bound}.

\paragraph{Proof of Lemma~\ref{lem:poten_h_bound}:} The function $h$ has the following derivatives
\begin{align*}
  \frac{d h}{d a} = \frac{1}{q} \parenth{a + 1}^{-1/q - 1}, \frac{d^2 h}{da^2} = -\frac{q+1}{q^2} \parenth{a + 1}^{-1/q - 2}.
\end{align*}
Using It\^o's formula, we obtain
\begin{align*}
  d h(\potentwo_t) &= \left.\frac{d h}{d a}\right|_{\potentwo_t} d\potentwo_t + \frac{1}{2} \left.\frac{d^2 h}{d a^2}\right|_{\potentwo_t} d\brackets{\potentwo}_t \\
  &= \frac{1}{q (\potentwo_t+1)^{1/q+1}} d\potentwo_t - \frac{1}{2} \frac{q+1}{q^2 (\potentwo_t+1)^{1/q+2}} \vecnorm{v_t}{2}^2 dt \\
  &\leq \frac{1}{q (\potentwo_t+1)^{1/q+1}} d\potentwo_t \\
  &\stackrel{(i)}{\leq} 64 q \alpha^2 \log(\dims) \dims^{2\beta-1/q} dt + \frac{v_t\tp dW_t}{q \parenth{\potentwo_t + 1}^{1/q+1}},
\end{align*}
where inequality (i) plugs in the bounds in Lemma~\ref{lem:potentwo_derivative_bound}.

Define a martingale $Y_t$ such that
\begin{align*}
  dY_t = \frac{v_t\tp dW_t}{q \parenth{\potentwo_t + 1}^{1/q+1}},
\end{align*}
with $Y_0 = 0$. According to the $\vecnorm{v_t}{2}$ upper bound in Lemma~\ref{lem:potentwo_derivative_bound}, we have
\begin{align*}
  \vecnorm{\frac{1}{q \parenth{\potentwo_t + 1}^{1 + 1/q}}v_t}{2}^2 &\leq 256.
\end{align*}
Hence the martingale $Y_t$ is well-defined.
According to the Dambis, Dubins-Schwarz theorem (see e.g. Revuz and Yor~\cite{revuz2013continuous} Section V.1 Theorem 1.7), there exits a Wiener process $\tilde{W}_t$ such that $Y_t$ has the same distribution as $\tilde{W}_{[Y]_t}$. Then we have for any $\gamma > 0$,
\begin{align}
  \label{eq:Y_t_tail_bound}
  \Prob \parenth{\max_{t \in [0, T]} Y_t \geq \gamma} \leq \Prob \parenth{\max_{t \in [0, T]} \tilde{W}_{256t} \geq \gamma} \leq \exp\parenth{-\frac{\gamma^2}{512 T}}.
\end{align}

Set $T = \frac{1}{32768 q \alpha^2 \log(\dims) d^{2\beta}}$ and $\Psi = \frac{1}{2} \parenth{\dims+1}^{-1/q}$. Observe that $\potentwo_0 = \dims$ and as a result $h(\potentwo_0) = -\parenth{\dims+1}^{-1/q}$. Then we have
\begin{align*}
  \Prob \parenth{\max_{t \in [0, T]}  h(\potentwo_t) \geq -\Psi}
  &\leq \Prob \parenth{\max_{t \in [0, T]}  Y_t \geq -\Psi + \parenth{\dims+1}^{-1/q} - \int_0^T 64q \alpha^2 \log(\dims) d^{2\beta - 1/q} dt }\\
  &\stackrel{(i)}{\leq} \Prob \parenth{\max_{t \in [0, T]}  Y_t \geq \frac{\Psi}{4}} \\
  &\stackrel{(ii)}{\leq} \exp\parenth{-\frac{\Psi^2}{8192T}} \\
  &\stackrel{(iii)}{\leq} \exp\parenth{-\frac{2}{3} q \alpha^2 \log(\dims) \dims^{2\beta - 2/q}} \notag \\
  &\stackrel{(iv)}{<} \frac{3}{10}.
\end{align*}
Inequality (i) follows from the choice of $T$. (ii) uses Equation~\eqref{eq:Y_t_tail_bound}. (iii) follows by plugging in $\Psi = \frac{1}{2}\parenth{\dims+1}^{-1/q}$ and $3^q \dims^2 \geq 2^q (\dims + 1)^2$. (iv) follows from $\beta q \geq 1$, $\dims \geq 3$, $q\geq 2$ and $3^{-4/3} < 0.3$.

\paragraph{Proof of Lemma~\ref{lem:poten_f_bound}:} The function $f$ has the following derivatives
\begin{align*}
  \frac{d f(a)}{d a} = \frac{1}{q} a^{1/q-1}, \frac{d^2 f(a, t)}{d a^2} = -\frac{q-1}{q^2} a^{1/q-2}.
\end{align*}
Using It\^o's formula, we obtain
\begin{align*}
  d f\parenth{\potentwo_t} &= \left.\frac{df}{da} \right|_{\potentwo_t} d\potentwo_t + \frac{1}{2} \left.\frac{d^2 f}{ d^2 a }\right|_{\potentwo_t} d \brackets{\potentwo}_t \\
  &= \frac{1}{q} \potentwo_t^{1/q-1} \parenth{v_t\tp dW_t + \delta_t dt} - \frac{q-1}{2q^2} \potentwo_t^{1/q-2} \vecnorm{v_t}{2}^2 dt.
\end{align*}
Using the bounds in Lemma~\ref{lem:potentwo_derivative_bound} and the martingale property of the term $\frac{1}{q} \potentwo_t^{1/q-1} v_t\tp dW_t$, we obtain
\begin{align*}
  d \Exs f(\potentwo_t) \leq \frac{2q}{t} \Exs f(\potentwo_t) dt.
\end{align*}
Solving the above differential equation, we obtain
\begin{align*}
  \Exs f(\potentwo_{t_2}) \leq \Exs f(\potentwo_{t_1}) \parenth{\frac{t_2}{t_1}}^{2q}, \forall t_2 > t_1 > 0.
\end{align*}
% subsection control_the_potential_in_two_time_periods (end)

\subsection{Proof of tensor bounds} % (fold)
\label{sub:proof_of_tensor_bounds}
In this subsection, we prove Lemma~\ref{lem:tensor_vector_bound},~\ref{lem:ttensor_A_isoperimetric_bound},~\ref{lem:ttensor_oneoverT_bound} and~\ref{lem:ttensor_swap}.
\paragraph{Proof of Lemma~\ref{lem:tensor_vector_bound}:}
Since $C$ is positive semi-definite, we can write its eigenvalue decomposition as follows $C = \sum_{i=1}^\dims \lambda_i v_i v_i\tp$, with $\lambda_i \geq 0$.
Then,
\begin{align*}
  &\vecnorm{\int B^{1/2} (x-\mu) (x-\mu)\tp C (x-\mu) p(x) dx}{2} \\
  &= \vecnorm{\sum_{i=1}^\dims \int B^{1/2} (x-\mu) \lambda_i \parenth{(x-\mu)\tp v_i}^2 p(x) dx}{2}\\
  &\stackrel{(i)}{\leq} \sum_{i=1}^\dims \lambda_i \vecnorm{\int B^{1/2} (x-\mu) \parenth{(x-\mu)\tp v_i}^2 p(x) dx}{2}\\
  &= \sum_{i=1}^\dims \lambda_i \max_{\vecnorm{\xi}{2}\leq 1} \int \xi\tp B^{1/2} (x-\mu) \parenth{(x-\mu)\tp v_i}^2 p(x) dx \\
  &\stackrel{(ii)}{\leq} \sum_{i=1}^\dims \lambda_i \max_{\vecnorm{\xi}{2}\leq 1} \parenth{\int \parenth{\xi\tp B^{1/2} (x-\mu)}^2 p(x) dx}^{1/2} \parenth{\int \parenth{(x-\mu)\tp v_i}^4 p(x) dx}^{1/2} \\
  &\stackrel{(iii)}{\leq} 16 \sum_{i=1}^\dims \lambda_i \max_{\vecnorm{\xi}{2}\leq 1} \parenth{\int \parenth{\xi\tp B^{1/2} (x-\mu)}^2 p(x) dx}^{1/2} \parenth{\int \parenth{(x-\mu)\tp v_i}^2 p(x) dx} \\
  &= 16\vecnorm{B^{1/2} A B^{1/2} }{2}^{1/2} \trace\parenth{A^{1/2}CA^{1/2}}
\end{align*}
Inequality (i) follows from triangular inequality. (ii) follows from Cauchy-Schwarz inequality. (iii) follows from the statement below, which upper bounds the fourth moment of a log-concave density via its second moment.

For any log-concave density $\nu$ and any vector $\theta\in \real^{\dims}$, we have
\begin{align}
  \label{eq:logconcave_regular}
  \parenth{\int \parenth{(x-\mu_\nu)\tp\theta}^a \nu(x) dx}^{1/a} \leq 2 \frac{a}{b} \parenth{\int \parenth{(x-\mu_\nu)\tp\theta}^b \nu(x) dx}^{1/b}
\end{align}
for $a \geq b > 0$, where $\mu_\nu$ is the mean of $\nu$. Equation~\eqref{eq:logconcave_regular} is proved e.g. in Corollary 5.7 of Gu\'edon et al.~\cite{guedon2014concentration} and the exact constant is provided in Proposition 3.8 of Lata{\l}a and Wojtaszczyk~\cite{latala2008infimum}.

In order to prove Lemma~\ref{lem:ttensor_A_isoperimetric_bound}, we need to introduce one additional lemma as follows.
\begin{lemma}
  \label{lem:trace_DeltaDelta_bound}
  Suppose that $\isoperi_k \geq \frac{1}{\alpha k^\beta}$ for all $k \leq \dims$ for some $0 < \beta \leq \frac{1}{2}$ and $\alpha \geq 1$.  For an isotropic log-concave density $p$ in $\real^\dims$ and a unit vector $v \in \real^\dims$, define $\Delta = \Exs_{X \sim p} \parenth{X\tp v} \cdot XX\tp$, then we have
  \begin{enumerate}
    \item For any orthogonal projection matrix $P \in \real^{\dims \times \dims}$ with rank $r$, we have
    \begin{align*}
      \trace\parenth{\Delta P \Delta} \leq 16 \isoperi^{-2}_{\min(2r, \dims)}
    \end{align*}
    \item For any positive semi-definite matrix $A$, we have
    \begin{align*}
      \trace\parenth{\Delta A \Delta} \leq 128 \alpha^2 \log(\dims) \parenth{\trace\parenth{A^{1/(2\beta)}}}^{2\beta}
    \end{align*}
  \end{enumerate}
\end{lemma}
This lemma was proved in Lemma 41 in an older version (arXiv version 2) of Lee and Vempala~\cite{lee2016eldan}. The main proof idea for the first part of Lemma~\ref{lem:trace_DeltaDelta_bound} appeared in Eldan~\cite{eldan2013thin} (Lemma 6). we provide a proof here for completeness.
\paragraph{Proof of Lemma~\ref{lem:trace_DeltaDelta_bound}:} For the first part, we have
\begin{align*}
  \trace\parenth{\Delta P \Delta} = \Exs_{X \sim p} X\tp \Delta P X \cdot X \tp v.
\end{align*}
Since $\Exs_{X\sim p} X\tp v = 0$, we can subtract the mean of the first term $X\tp \Delta P X$ without changing the value of $\trace\parenth{\Delta P \Delta} $. Then
\begin{align*}
  \trace\parenth{\Delta P \Delta} &= \Exs_{X\sim p} \brackets{\parenth{X\tp \Delta P X  - \Exs_{Y \sim p} Y \tp \Delta P Y} \cdot X \tp v} \\
  &\stackrel{(i)}{\leq} \parenth{\Exs_{X\sim p}(X\tp v)^2}^{1/2} \parenth{\Var_{X \sim p }\parenth{X \tp \Delta P X}}^{1/2} \\
  &\stackrel{(ii)}{\leq} 2 \isoperi_{\min(2r, \dims)}^{-1} \parenth{\Exs_{X \sim p} \vecnorm{\Delta P X + P\tp \Delta\tp X}{2}^2}^{1/2} \\
  & \leq 4 \isoperi_{\min(2r, \dims)}^{-1} \parenth{\trace\parenth{\Delta P \Delta}}^{1/2}.
\end{align*}
Inequality (i) follows from the Cauchy-Schwarz inequality. Inequality (ii) follows from the fact that $\Exs_{X\sim p}(X\tp v)^2 = 1$ as $p$ is isotropic and that the inverse Poincar\'e constant is upper bounded by twice of inverse of the squared isoperimetric coefficient (also known as Cheeger's inequality~\cite{maz1960classes,cheeger1969lower} or Theorem 1.1 in Milman~\cite{milman2009role}). The matrix $\Delta P + P\tp \Delta$ has rank at most $\min(2r, \dims)$. Rearranging the terms in the above equation, we conclude the first part of Lemma~\ref{lem:trace_DeltaDelta_bound}.

For the second part, we write the matrix $A$ in its eigenvalue decomposition and group the terms by eigenvalues. We have
\begin{align*}
  A = \sum_{i=1}^\dims \lambda_i v_i v_i\tp = \sum_{j=1}^J A_j + B,
\end{align*}
where $A_i$ has eigenvalues between the interval $(\vecnorm{A}{2} e^{i-1} /\dims, \vecnorm{A}{2} e^{i} /\dims]$ and $B$ has eigenvalues smaller than or equal to $\vecnorm{A}{2}/\dims$. Because the intervals have right bounds increasing exponentially, we have $J = \lceil \log(\dims) \rceil$. Let $P_i$ be the orthogonal projection matrix formed by the eigenvectors in $A_i$. Then we have
\begin{align}
  \label{eq:trace_DeltaDelta_Ai_bound}
  \trace\parenth{\Delta A_i \Delta} \leq \vecnorm{A_i}{2} \trace\parenth{\Delta P_i \Delta} \stackrel{(i)}{\leq} 16 \vecnorm{A_i}{2} \isoperi^{-2}_{\min(2 \text{rank}(A_i), \dims)} \stackrel{(ii)}{\leq} 16 \alpha^2 \vecnorm{A_i}{2} \cdot \parenth{2 \text{rank}(A_i)}^{2\beta},
\end{align}
where inequality (i) follows from the first part of Lemma~\ref{lem:trace_DeltaDelta_bound} and inequality (ii) follows from the hypothesis of Lemma~\ref{lem:trace_DeltaDelta_bound}.
Similarly for matrix $B$, we have
\begin{align}
  \label{eq:trace_DeltaDelta_B_bound}
  \trace\parenth{\Delta B \Delta} \stackrel{(i)}{\leq} 16 \alpha^2 \vecnorm{B}{2} \parenth{2\text{rank}(B)}^{2\beta} \stackrel{(ii)}{\leq} 32 \alpha^2 \vecnorm{A}{2},
\end{align}
where inequality (i) follows from the hypothesis of Lemma~\ref{lem:trace_DeltaDelta_bound} and inequality (ii) follows from the fact that $\vecnorm{B}{2} \leq \vecnorm{A}{2}/d$ and $2\beta \leq 1$.
Putting the bounds~\eqref{eq:trace_DeltaDelta_Ai_bound} and~\eqref{eq:trace_DeltaDelta_B_bound} together, we have
\begin{align*}
  \trace\parenth{\Delta A \Delta} &= \sum_{j=1}^J \trace\parenth{\Delta A_j\Delta} + \trace\parenth{\Delta  B \Delta} \\
  &\leq 16 \alpha^2 \parenth{\sum_{j=1}^J \vecnorm{A_j}{2} \cdot \parenth{2\text{rank}(A_j)}^{2\beta} + 2\vecnorm{A}{2}} \\
  &\stackrel{(i)}{\leq} 16 \alpha^2 \brackets{\parenth{\sum_{j=1}^J \vecnorm{A_j}{2}^{1/(2\beta)} \cdot \parenth{2\text{rank}(A_j)}}^{2\beta} \cdot \parenth{J}^{1-2\beta} +  2\vecnorm{A}{2}}\\
  &\stackrel{(ii)}{\leq} 16 \alpha^2\brackets{\parenth{2 e \trace\parenth{A^{1/(2\beta)}}}^{2\beta} \cdot \parenth{J}^{1-2\beta} + 2 \vecnorm{A}{2}} \\
  &\leq 128 \alpha^2 \log(\dims) \parenth{\trace\parenth{A^{1/(2\beta)}}}^{2\beta}.
\end{align*}
Inequality (i) follows from Holder's inequality and inequality (ii) follows from the fact that $\vecnorm{A_j}{2}^{1/2\beta} \text{rank}(A_j) \leq e \trace\parenth{A_{j}^{1/2\beta}}$ due to the construction of $A_j$. This concludes the second part of Lemma~\ref{lem:trace_DeltaDelta_bound}.

\paragraph{Proof of Lemma~\ref{lem:ttensor_A_isoperimetric_bound}:} Let $\mu$ be the mean of $p$. First, for $X$ a random vector in $\real^\dims$ drawn from $p$, we define the standardized random variable $A^{-1/2} (X - \mu)$ and its density $\denb$. $\denb$ is an isotropic log-concave density. Then through a change of variable, we have
\begin{align*}
  &\ttensor_p \parenth{A^{q-2}, \Ind_\dims, \Ind_\dims} \\
  =& \int\int (x-\mu)\tp A^{q-2} (y-\mu) \cdot (x-\mu)\tp (y-\mu) \cdot (x-\mu) \tp (y-\mu) p(x) p(y) dx dy \\
  =& \int\int \parenth{x\tp A^{q-1} y} (x\tp A y) (x \tp A y) \denb(x) \denb(y) dx dy \\
  \leq& \int\int \parenth{x\tp A^{q} y} (x\tp A y) (x \tp y) \denb(x) \denb(y) dx dy \\
  =& \ttensor_\denb \parenth{A^{q}, A, \Ind_\dims},
\end{align*}
where the last inequality follows from Lemma~\ref{lem:ttensor_swap}. $A^q$ is positive semi-definite and we write down its eigenvalue decomposition $A^q =  \sum_{i=1}^\dims \lambda_i v_i v_i \tp$ with $\lambda_i \geq 0$. Since $\denb$ is isotropic, we can rewrite the 3-Tensor into a summation form and apply Lemma~\ref{lem:trace_DeltaDelta_bound}.
\begin{align*}
  &\ttensor_\denb \parenth{A^{q}, A, \Ind_\dims} \\
  & = \int \int \parenth{x \tp A^q y} \parenth{x \tp A y} \parenth{x\tp y} \denb(x) \denb(y) dx dy  \\
  & = \sum_{i=1}^\dims \lambda_i \int \int \parenth{x \tp v_i} \parenth{y \tp v_i} \parenth{x \tp A y} \parenth{x\tp y} \denb(x) \denb(y) dx dy  \\
  & = \sum_{i=1}^\dims \lambda_i \trace\parenth{\Delta_i A \Delta_i} \\
  &\stackrel{(i)}{\leq}  128 \alpha^2 \log(\dims) \parenth{\trace(A^{1/2\beta})}^{2\beta} \parenth{\sum_{i=1}^\dims \lambda_i} \\
  &= 128 \alpha^2 \log(\dims) \parenth{\trace(A^{1/2\beta})}^{2\beta} \trace(A^q) \\
  &\stackrel{(ii)}{\leq} 128 \alpha^2 \log(\dims) \trace(A^q) \brackets{\trace\parenth{A^q}^{1/(2\beta q)} \parenth{\dims}^{1 - 1/(2\beta q)}}^{2\beta} \\
  &= 128 \alpha^2 \log(\dims) \dims^{2\beta - 1/q} \trace(A^q) ^{1 + 1/q},
\end{align*}
where we define $\Delta_i = \int (x\tp v_i) x x\tp \varrho(x) dx$, inequality (i) follows from Lemma~\ref{lem:trace_DeltaDelta_bound} and that $\varrho$ is isotropic, inequality (ii) follows from Cauchy-Schwarz inequality and the assumption that $q \geq \frac{1}{2\beta}$.

\paragraph{Proof of Lemma~\ref{lem:ttensor_oneoverT_bound}:}
Without loss of generality, we can assume that the density $p$ has mean $0$. Its covariance matrix $A$ is positive semi-definite and invertible. We can write down its eigenvalue decomposition as follows $A = \sum_{i=1}^\dims \lambda_i v_i v_i\tp$ with $\lambda_i > 0$ and $v_i$ are eigenvectors with norm $1$. Then $A^{q}$ has an eigenvalue decomposition with the same eigenvectors $A^q = \sum_{i=1}^\dims \lambda_i^q v_i v_i\tp$.
Define $\Delta_i = \Exs_{X \sim p} (X\tp A^{-1/2}v_i) X X \tp$, then
\begin{align}
  \label{eq:ttensor_oneoverT_first_step}
  \ttensor_p\parenth{A^{q-2}, \Ind_\dims, \Ind_\dims} &= \Exs_{X, Y \sim p} \parenth{X\tp A^{q-2} Y} (X\tp Y) (X \tp Y) \notag \\
  &= \sum_{i=1}^\dims \lambda_i^{q-1} \trace\parenth{\Delta_i \Delta_i}.
\end{align}
Next we bound the terms $\trace\parenth{\Delta_i \Delta_i}$. We have
\begin{align*}
  \trace\parenth{\Delta_i \Delta_i} &= \Exs_{X \sim p} \parenth{X \tp A^{-1/2} v_i} X\tp \Delta_i X \\
  &\stackrel{(i)}{=}  \Exs_{X \sim p} \parenth{X \tp A^{-1/2} v_i} \parenth{X\tp \Delta_i X - \Exs_{Y \sim p} \brackets{Y\tp \Delta_i Y}} \\
  &\stackrel{(ii)}{\leq} \parenth{\Exs_{X \sim p} \parenth{X \tp A^{-1/2} v_i}^2 }^{1/2} \parenth{\Var\parenth{X \tp \Delta_i X}}^{1/2} \\
  & \stackrel{(iii)}{=} \parenth{\Var_{X \sim p}\parenth{X \tp \Delta_i X}}^{1/2} \\
  &\stackrel{(iv)}{\leq} \parenth{\Exs_{X \sim p} \frac{1}{\tau} \vecnorm{\Delta_i X + \Delta_i X}{2}^2}^{1/2} \\
  &\stackrel{(v)}{\leq} \parenth{\frac{4}{\tau} \trace\parenth{A \Delta_i \Delta_i}}^{1/2}.
\end{align*}
Equality (i) is because $\Exs_{X \sim p} X = 0$. Inequality (ii) follows from Cauchy-Schwarz inequality. Equality (iii) follows from the definition of the covariance matrix $\Exs_{X\sim p} XX\tp = A$. Inequality (iv) follows from the Brascamp-Lieb inequality (or Hessian Poincar\'e, see Theorem 4.1 in Brascamp and Lieb~\cite{brascamp2002extensions}) together with the assumption that $p$ is more log-concave than $\Normal(0, \frac{1}{\tau}\Ind_\dims)$.

Plugging the bounds of the terms $\trace\parenth{\Delta_i \Delta_i}$ into Equation~\eqref{eq:ttensor_oneoverT_first_step}, we obtain
\begin{align*}
  \ttensor_p\parenth{A^{q-2}, \Ind_\dims, \Ind_\dims} &= \sum_{i=1}^\dims \lambda_i^{q-1} \trace\parenth{\Delta_i \Delta_i} \\
  &\leq \sum_{i=1}^\dims \lambda_i^{q-1}  \parenth{\frac{4}{\tau} \trace\parenth{A \Delta_i \Delta_i}}^{1/2} \\
  &\stackrel{(i)}{\leq} \frac{2}{\tau^{1/2}} \parenth{\sum_{i=1}^\dims \lambda_i^{q}}^{1/2}  \parenth{\sum_{i=1}^\dims \lambda_i^{q-2} \trace\parenth{A \Delta_i \Delta_i} }^{1/2} \\
  &= \frac{2}{\tau^{1/2}} \parenth{\trace\parenth{A^q}}^{1/2}  \parenth{\Exs_{X, Y \sim p} \parenth{X\tp A^{q-3} Y} (X\tp A Y) (X \tp Y) }^{1/2} \\
  &\stackrel{(ii)}{\leq} \frac{2}{\tau^{1/2}} \parenth{\trace\parenth{A^q}}^{1/2}  \parenth{\Exs_{X, Y \sim p} \parenth{X\tp A^{q-2} Y} (X\tp Y) (X \tp Y) }^{1/2} \\
  &= \frac{2}{\tau^{1/2}} \parenth{\trace\parenth{A^q}}^{1/2}  \brackets{\ttensor_p\parenth{A^{q-2}, \Ind_\dims, \Ind_\dims} }^{1/2}.
\end{align*}
Inequality (i) follows from Cauchy-Schwarz inequality. For $q \geq 3$, inequality (ii) follows from Lemma~\ref{lem:ttensor_swap}. From the above equation, after rearranging the terms, we obtain
\begin{align*}
  \ttensor_p\parenth{A^{q-2}, \Ind_\dims, \Ind_\dims} \leq \frac{4}{\tau} \trace\parenth{A^q}.
\end{align*}

\paragraph{Proof of Lemma~\ref{lem:ttensor_swap}: } This lemma is proved in Lemma 43 in an older version (arXiv version 2) of Lee and Vempala~\cite{lee2016eldan}, we provide a proof here for completeness.

Without loss of generality, we can assume that the density $p$ has mean $0$. For $i \in \braces{1, \cdots, \dims}$, we define $\Delta_i = \Exs_{X\sim p} B^{1/2} X X \tp B^{1/2} X\tp C^{1/2} e_i$ where $e_i \in \real^\dims$ is the vector with $i$-th coordinate 1 and 0 elsewhere. We have $\sum_{i=1}^\dims e_i e_i \tp  = \Ind_\dims$. We can rewrite the tensor on the left hand side as a sum of traces.
\begin{align}
  \label{eq:proof_lemma_swap_main}
  &\ttensor_{p}(B^{1/2}A^\delta B^{1/2}, B^{1/2}A^{1-\delta}B^{1/2}, C) \notag \\
  =& \Exs_{X, Y \sim p} X\tp B^{1/2}A^\delta B^{1/2} Y \cdot X\tp B^{1/2}A^{1-\delta}B^{1/2} Y \cdot X \tp C Y \notag \\
  =& \sum_{i=1}^\dims \Exs_{X, Y \sim p} X\tp B^{1/2}A^\delta B^{1/2} Y \cdot X\tp B^{1/2}A^{1-\delta}B^{1/2} Y \cdot X\tp C^{1/2} e_i \cdot Y\tp C^{1/2} e_i \notag \\
  =& \sum_{i=1}^\dims \trace\parenth{A^{\delta} \Delta_i A^{1-\delta} \Delta_i}.
\end{align}
For any symmetric matrix $F$, a positive-semidefinite matrix $G$ and $\delta \in [0, 1]$, we have
\begin{align}
  \label{eq:trace_ineq_lieb_thirring}
  \trace\parenth{G^\delta F G^{1-\delta} F} \leq \trace\parenth{G F^2}.
\end{align}
Applying the above trace inequality~\eqref{eq:trace_ineq_lieb_thirring} that we prove later for completeness (see also Lemma 2.1 in Zhu et al.~\cite{allen2016using}),
we obtain
\begin{align*}
  \trace\parenth{A^{\delta} \Delta_i A^{1-\delta} \Delta_i} \leq \trace\parenth{A \Delta_i \Delta_i}.
\end{align*}
Writing the sum of traces in Equation~\eqref{eq:proof_lemma_swap_main} back to the 3-Tensor form, we conclude Lemma~\ref{lem:ttensor_swap}.

It remains to prove the trace inequality in Equation~\eqref{eq:trace_ineq_lieb_thirring}. Without loss of generality, we can assume $G$ is diagonal. Hence, we have
\begin{align*}
  \trace\parenth{G^\delta F G^{1-\delta} F} &= \sum_{i = 1}^\dims \sum_{j = 1}^\dims G_{ii}^\delta G_{jj}^{1-\delta} F_{ij}^2 \\
  &\leq \sum_{i=1}^\dims \sum_{j=1}^\dims \parenth{\delta G_{ii} + (1-\delta) G_{jj}} F_{ij}^2 \\
  &= \delta \sum_{i=1}^\dims \sum_{j=1}^\dims G_{ii} F_{ij}^2 + (1-\delta) \sum_{i=1}^\dims \sum_{j=1}^\dims G_{jj} F_{ij}^2 \\
  &= \trace\parenth{G F^2},
\end{align*}
where the inequality follows from Jensen's inequality and the fact that the logarithm function is concave (or the inequality of arithmetic and geometric means).
% subsection proof_of_tensor_bounds (end)

% section main_results (end)

% \section{Discussion} % (fold)
% \label{sec:discussion}

% % section discussion (end)

%%%%%%%%%%%%%%%%%%%%%%%%%%%%%%%%%%%%%%%%%%%%%%%%%%%%%%%%%%%%%%%%%%%%%%
\subsection*{Acknowledgements}
Yuansi Chen has received funding from the European Research Council under
the Grant Agreement No 786461 (CausalStats - ERC-2017-ADG). We acknowledge
scientific interaction and exchange at ``ETH Foundations of Data Science''. We thank Peter B\"uhlmann and Bin Yu for their continuous support and encouragement. We thank Afonso Bandeira, Raaz Dwivedi, Ronen Eldan, Yin Tat Lee and Martin Wainwright for helpful discussions. We thank Bo'az Klartag and Joseph Lehec for pointing out a mistake in the previous revision. We also thank anonymous reviewers for their careful reading of our manuscript and their suggestions on presentation and writing.
%% APPENDIX %%%%%%%%%%%%%%%%%%%%%%%%%%%%%%%%%%%%%%%%%%%%%%%%%%%%%%%%%%%%

\appendix
%!TEX root = main_paper.tex

\section{Proof of Lemma~\ref{lem:existence_SDE_solution} and derivatives} % (fold)
\label{sec:proof_of_derivatives}
In this section, we first prove the existence and uniqueness of the SDE solution in Lemma~\ref{lem:existence_SDE_solution} and then derive the derivatives of $p_t$, $A_t$ and $\potentwo_t$ in Equation~\eqref{eq:derivative_pt}, Equation~\eqref{eq:known_SDE_change_At} and~\eqref{eq:poten2_derivative} using It\^o's calculus. Similar results are also proved in Eldan~\cite{eldan2013thin} and Lee and Vempala~\cite{lee2016eldan} since a similar stochastic localization is used. We provide a proof here for completeness.

\paragraph{Proof of Lemma~\ref{lem:existence_SDE_solution}:}
We can rewrite the stochastic differential equation~\eqref{eq:def_SDE_ctBt} as follows to make the dependency clear:
\begin{align*}
  dc_t &= \co^{-1/2} dW_t + \co^{-1} \mu\parenth{c_t, B_t} dt\\
  dB_t &= \co^{-1} dt,
\end{align*}
where
\begin{align*}
  \mu(c, B) &= \int x \denb(c, B, x) dx, \\
  \denb(c, B, x) &= \frac{e^{c\tp x - \frac{1}{2}x\tp B x} p(x)}{\int_{\real^\dims} e^{c \tp x - \frac{1}{2}y\tp B y} p(y) dy}.
\end{align*}
Since $p$ has a compact support, given $x \in \real^\dims$, $\denb(\cdot, \cdot, x)$ as a function of $(c, B)$ is Lipschitz in $c$ and $B$. Similarly, $\mu$ is also Lipschitz in $c$ and $B$. Consequently, $\co^{-1/2}$, $\co \mu(c_t, B_t)$ and $\co^{-1}$ are all bounded and Lipschitz on $c_t$ and $B_t$ on the compact support. Applying the existence and uniqueness theorem of SDE solutions (Theorem 5.2 in {\O}ksendal~\cite{oksendal2003stochastic}), we show that the SDE solution exists and is unique on the time interval $[0, T]$ for any $T > 0$.

Next, we derive the derivative of $p_t$. Define
\begin{align*}
  G_t(x) &= e^{c_t\tp x - \frac{1}{2} x\tp B_t x} p(x), \\
  V_t & = \int G_t(x) dx.
\end{align*}
Then $p_t(x)$ can be written as $\frac{G_t(x)}{V_t}$. Let $S_t(x)$ denote the quadratic variation of the process $c_t\tp x$. We have
\begin{align*}
  d S_t(x) = x\tp \co^{-1} x dt.
\end{align*}
Using It\^o's formula, we have
\begin{align*}
  dG_t(x) &= \parenth{x\tp (dc_t) - \frac{1}{2} x\tp dB_t x + \frac{1}{2} dS_t} G_t(x) \\
  &= \parenth{x\tp \co^{-1/2}dW_t + x\tp \co^{-1} \mu_t dt} G_t(x), \\
  dV_t &= \int dG_t(x) dx = V_t \parenth{\mu_t\tp \co^{-1/2}dW_t + \mu_t\tp \co^{-1} \mu_t dt}.
\end{align*}
Using It\^o's formula on the inverse of $V_t$, we have
\begin{align*}
  d V_t^{-1} &= -\frac{dV_t}{V_t^2} + \frac{d \brackets{V}_t} {V_t^3} \\
  &= - V_t^{-1} \brackets{\mu_t\tp \co^{-1/2} dW_t + \mu_t\tp \co^{-1} \mu_t dt} + V_t^{-1} \mu_t\tp \co^{-1} \mu_t dt \\
  &= - V_t^{-1} \mu_t\tp \co^{-1/2} dW_t.
\end{align*}
Using It\^o's formula on $p_t$, with the above derivatives, we obtain
\begin{align*}
  dp_t(x) &= d \parenth{V_t^{-1} G_t(x)} \\
          &= \parenth{G_t(x) dV_t^{-1} + V_t^{-1}dG_t(x) + d\brackets{V^{-1}, G(x)}_t} \\
          &= \parenth{x-\mu_t}\tp \co^{-1/2} dW_t p_t(x).
\end{align*}

Then we derive the derivative of $A_t$. By the definition of $A_t$, we have
\begin{align*}
  A_t = \int \parenth{x - \mu_t} \parenth{x - \mu_t}\tp p_t(x) dx,
\end{align*}
where $\mu_t = \int_{\real^\dims} x p_t(x) dx$. Using It\^o's formula on $\mu_t$, we obtain
\begin{align*}
  d\mu_t &= \int x d p_t(x) dx \\
  &= \int x (x- \mu_t)\tp \co^{-1/2} dW_t p_t(x) dx \\
  &= \int (x - \mu_t) (x- \mu_t)\tp \co^{-1/2} dW_t dx\\
  &= A_t \co^{-1/2} dW_t.
\end{align*}
Using It\^o's formula on $A_t$ and viewing it as a function of $\mu_t$ and $p_t$, we obtain
\begin{align*}
  dA_t =& \int \parenth{x - \mu_t}\parenth{x - \mu_t}\tp dp_t(x) dx - \int d\mu_t\parenth{x - \mu_t}\tp p_t(x) dx - \int \parenth{x - \mu_t}\parenth{d\mu_t}\tp p_t(x) dx \\
  &-\frac{1}{2}\cdot 2 \int \parenth{x - \mu_t} d\brackets{\mu_t\tp, p_t(x)}_t dx - \frac{1}{2}\cdot 2 \int  d\brackets{\mu_t, p_t(x)}_t \parenth{x - \mu_t}\tp  dx \\
  &+ \frac{1}{2} \cdot 2 d\brackets{\mu_t, \mu_t\tp}_t \int p_t(x) dx.
\end{align*}
We observe that $\int d\mu_t\parenth{x - \mu_t}\tp p_t(x) dx = 0$ and $\int \parenth{x - \mu_t}\parenth{d\mu_t}\tp p_t(x) dx = 0$. Then,
\begin{align*}
  d\brackets{\mu_t\tp, p_t(x)}_t &= \parenth{x - \mu_t}\tp \co^{-1} A_t p_t(x) dt,\\
  d\brackets{\mu_t, p_t(x)}_t &= A_t  \co^{-1} \parenth{x - \mu_t} p_t(x) dt, \\
  d\brackets{\mu_t, \mu_t}_t &= A_t A^{-1} A_t dt.
\end{align*}
Combining all the terms together, we have
\begin{align*}
  dA_t = \int \parenth{x - \mu_t} \parenth{x - \mu_t}\tp \parenth{\parenth{x-\mu_t}\tp \co^{-1/2} dW_t} p_t(x) dx - A_t\co^{-1} A_t dt.
\end{align*}

Finally, we derive the derivative of $\potentwo_t$. Define the function $\Gamma: \real^{\dims \times \dims} \mapsto \real$ as $\Gamma(X) = \trace\parenth{X^q}$. The first-order and second-order derivatives of $\Gamma$ are given by
\begin{align*}
  \left.\frac{\partial \Gamma}{\partial X}\right|_{H} = q \trace\parenth{ X^{q-1} H}, \left.\frac{\partial^2 \Gamma}{\partial X \partial X}\right|_{H_1, H_2} = q \sum_{a=0}^{q-2} \trace\parenth{X^a H_2 X^{q-2-a} H_1}.
\end{align*}
Using the above derivatives and It\^o's formula, we obtain
\begin{align}
  \label{eq:poten_two_derivative_intermediate}
  d\potentwo_t = d \trace\parenth{\qo_t^q} = q \trace\parenth{ \qo_t^{q-1} d \qo_t} + \frac{q}{2} \sum_{a = 0}^{q-2} \sum_{i,j,k,l=1}^{\dims} \trace\parenth{\qo_t^a E_{ij} \qo_t^{q-2-a} E_{kl}} d\brackets{Q_{ij}, Q_{kl}}_t,
\end{align}
where $E_{ij}$ is the matrix that takes $1$ at the entry $(i, j)$ and $0$ otherwise and $Q_{ij, t}$ is the stochastic process defined by the $(i,j)$ entry of $Q_t$.
Using the derivative of $A_t$ in Equation~\eqref{eq:known_SDE_change_At}, we have
\begin{align*}
  d\qo_t &= \int  \co^{-1/2} \parenth{x - \mu_t} \parenth{x - \mu_t}\tp \co^{-1/2} \parenth{\parenth{x-\mu_t}\tp \co^{-1/2} dW_t} p_t(x) dx - \qo_t^2 dt, \\
  d\brackets{\qo_{ij}, \qo_{kl}}_t &= \int \int z(x)_i z(x)_j z(y)_k z(y)_l (x-\mu_t)\tp \co^{-1} (y-\mu_t) p_t(x) p_t(y)dx dy dt,
\end{align*}
where $z(x)_i$ is the $i$-th coordinate of $\brackets{\co^{-1/2}(x-\mu_t)}$.
Plugging the expressions of $dA_t$ and $d\brackets{A_{ij}, A_{kl}}_t $ into Equation~\eqref{eq:poten_two_derivative_intermediate}, we obtain
\begin{align*}
  d\potentwo_t &=  q \int \parenth{x-\mu_t}\tp \co^{-1/2} \parenth{\qo_t}^{q-1} \co^{-1/2} \parenth{x-\mu_t} \parenth{x-\mu_t}\tp \co^{-1/2} dW_t p_t(x) dx \notag \\
  &- q \trace\parenth{\qo_t^{q+1} } dt \notag \\
  &+ \frac{q}{2} \sum_{a = 0}^{q-2} \int \int \parenth{x-\mu_t}\tp \co^{-1/2} \qo_t^{a} \co^{-1/2} \parenth{y-\mu_t} \notag \\
  &\cdot \parenth{x-\mu_t}\tp \co^{-1/2} \qo_t^{q-2-a} \co^{-1/2} \parenth{y-\mu_t} \parenth{x-\mu_t}\tp \co^{-1} \parenth{y-\mu_t} p_t(x) p_t(y) dx dy dt.
\end{align*}

% section proof_of_derivatives (end)

%%%%%%%%%%%% DO BIBLIOGRAPHY %%%%%%%%%%%%%%%%%%%%%%%%%%%%%%%%%%%%%%%%%%

\bibliographystyle{abbrv}
\bibliography{ref}

%%%%%%%%%%%%%%%%%%%%%%%%%%%%%%%%%%%%%%%%%%%%%%%%%%%%%%%%%%%%%%%%%%%%%%%%
\end{document}